\documentclass[11pt]{article}                     
\usepackage[margin = 1in]{geometry}
\usepackage{multicol}
\usepackage{amsmath,amssymb,latexsym,amsfonts,epsfig, color,geometry,tabularx,float,hyperref, xfrac, bm, booktabs, amsthm, mathrsfs}
\usepackage{tikz}
\usepackage{enumitem}
\usepackage{multirow}
\usepackage[numbers,sort]{natbib}
\newcommand\R{{\mathbb R}}
\newcommand\N{{\mathbb N}}
\newcommand\cK{{\mathcal{K}}}
\newcommand\cM{{\mathbb{M}}}
\newcommand\cP{{\mathbb{P}}}

\newcommand\SOMS{{\operatorname{SOMS}}}
\newcommand\SOS{{\operatorname{SOS}}}
\newcommand\SONC{{\operatorname{SONC}}}
\newcommand\SAGEmath{{\operatorname{SAG}}}
\newcommand\SAGEtext{SAG}
\newcommand\SAGElongtext{sums of arithmetic-mean/geometric-mean}
\newcommand\AGElongtext{arithmetic-mean/geometric-mean}
\newcommand\AGEtext{AM/GM}
\newcommand\DSOS{{\operatorname{DSOS}}}
\newcommand\SDSOS{{\operatorname{SDSOS}}}

\newcommand\KSOMS{\cK_\SOMS}
\newcommand\KSOS{\cK_\SOS}
\newcommand\KSONC{\cK_\SONC}
\newcommand\KSAGE{\cK_\SAGEmath{}}
\newcommand\KDSOS{\cK_\DSOS}
\newcommand\KSDSOS{\cK_\SDSOS}

\newcommand{\germans}{Positivstellens\"atze}
\newcommand{\german}{Positivstellensatz}
\newcommand{\schmu}{Schm\"{u}dgen}

\newcommand{\tr}{^\intercal}

\newcommand{\B}{{\mathcal B}}

\newtheorem{myexample}{Example}
\newtheorem{myremark}{Remark}
\newtheorem{theorem}{Theorem}
\newtheorem{proposition}{Proposition}
\newtheorem{corollary}{Corollary}
\newtheorem{definition}{Definition}
\newtheorem{lemma}{Lemma}

\title{Sparse non-SOS Putinar-type Positivstellens\"atze}

\author{ \small
Lorenz M. Roebers\thanks{\scriptsize Tilburg School of Economics and Management, Tilburg University, The Netherlands, {\tt l.m.roebers@tilburguniversity.edu}} 
\and Juan C. Vera\thanks{\scriptsize Tilburg School of Economics and Management, Tilburg University, The Netherlands, {\tt j.c.veralizcano@uvt.nl}}
\and Luis F. Zuluaga\thanks{\scriptsize Department of Industrial and Systems Engineering, Lehigh University, USA, {\tt
luis.zuluaga@lehigh.edu}}}

\begin{document}


\maketitle

\begin{abstract}

Recently, non-SOS Positivstellens\"atze for polynomials on compact semialgebraic sets, following the general form of Schm\"{u}dgen's Positivstellensatz, have been derived by appropriately replacing the SOS polynomials with other classes of polynomials. An open question in the literature is how to obtain similar results following the general form of Putinar's Positivstellensatz. Extrapolating the algebraic geometry tools used to obtain this type of result in the SOS case fails to answer this question, because algebraic geometry methods strongly use hallmark properties of the class of SOS polynomials, such as closure under multiplication and closure under composition with other polynomials.
In this article, using a new approach, we show the existence of Putinar-type  Positivstellens\"atze that are constructed using non-SOS classes of non-negative polynomials, such as SONC, SDSOS and DSOS polynomials. Even not necessarily non-negative classes of polynomials such as \SAGElongtext{} polynomials could be used.  Furthermore, we show that these certificates can be written with inherent sparsity characteristics. Such characteristics can be further exploited when the sparsity structure of both the polynomial whose non-negativity is being certified and the polynomials defining the semialgebraic set of interest are known. In contrast with related literature focused on exploiting sparsity in SOS Positivstellens\"atze, these latter results show how to exploit sparsity in a more general setting in which non-SOS polynomials are used to construct the Positivstellens\"atze.
\end{abstract}


\section{Introduction}

Let us refer to Positivstellens\"atze (i.e., certificates of nonegativity) such as Schm\"{u}dgen's~\cite{schmudgen1991k}, Putinar's~\cite{putinar1993positive}, and Reznick's~\citep{reznick1995uniform} Positivstellensatz, in which sum of squares (SOS) polynomials are used as the base class of polynomials to certify non-negativity, as {\em SOS Positivstellens\"atze}. This type of \germans{} have been studied for more than a century in algebraic geometry~\citep[see, e.g.][]{reznick1995uniform};
 together with the dual theory of moments, they are the foundation of modern {\em polynomial optimization}~\citep[see, e.g.,][]{lasserre2009moments, anjo12}. Broadly speaking, a Positivstellensatz can be translated into an approximation hierarchy for a polynomial optimization problem via a restricted-degree version of the Positivstellensatz~\citep[see, e.g.,][]{lasserre2009moments}. In turn, the ability to solve the resulting optimization problems in the approximation hierarchy depends on different characteristics of the Positivstellensatz, including number of terms, sparsity structure, and the base class of polynomials used to certify non-negativity (c.f., SOS polynomials in SOS Positivstellens\"atze).
A theorem guaranteeing the existence of a Positivstellensatz implies the convergence of its associated approximation hierarchy~\citep[see, e.g.,][]{lasserre2009moments}; explicit bounds on the degree of the polynomials required for such Positivstellensatz implies a rate of convergence for the hierarchy bounds~\citep[see, e.g.,][]{KlLaur10-degreeBounds}.

In what follows, we answer, in more generality, an important open question in the literature regarding \germans{}~\citep[][Sec.~6]{dressler2017positivstellensatz}. Namely, we prove the existence of a new class of {\em non-SOS \germans{}} that resembles Putinar's \german{}, but do not require the use of SOS polynomials to certify non-negativity (see Theorem~\ref{thm:semiSparse}). In particular, well known classes of non-SOS polynomials can be used to certify the non-negativity in the \germans{} (see Corollary~\ref{cor:semiSparse}). Moreover, much like, and going beyond results that aim at exploiting sparsity in SOS \germans{}~\citep[see, e.g.,][]{mai2020sparse, lasserre2006convergent}, we show that sparsity can also be exploited in the proposed non-SOS \germans{} (see Theorem~\ref{thm:mainSparse}), and in particular, when considering well known classes of non-SOS polynomials (see Corollary~\ref{cor:mainSparse}). Next, we discuss these results in more detail before formally stating and proving them thereafter.

Let $g_j$, $j=1,\dots,m$, be given polynomials such that the basic semialgebraic set $S = \{ x \in \R^n: g_1(x) \ge 0, \ldots, g_m(x) \ge 0\}$ is compact.
A fundamental SOS Positivstellensatz is Schm\"{u}dgen's Positivstellensatz~\cite{schmudgen1991k}, which states that every polynomial $p$ that is positive on S,
can be written in the form
\begin{equation}
\label{eq:schmudgenform}
p(x) = \sum_{\alpha \in \{0,1\}^m} \sigma_\alpha(x) g_1^{\alpha_1}(x) \cdots  g_m^{\alpha_m}(x),
\end{equation}
where $\sigma_\alpha$, for all $\alpha \in \{0,1\}^m$, are SOS polynomials.
As SOS polynomials are non-negative everywhere and $g_j(x) \ge 0$ for all $x \in S$, $j=1,\dots,m$, the right hand side of~\eqref{eq:schmudgenform} makes the non-negativity of $p$ on~$S$ evident.

Membership in the class of SOS polynomials, even though semidefinite programming (SDP) representable~\citep[see, e.g.,][]{BlekPT13}, is computationally expensive in practical terms~\citep[see, e.g.,][]{lasserre2006convergent}.
This fact has motivated recent research to focus on obtaining non-SOS Positivstellens\"atze in which non-SOS polynomials are used as the base class of polynomials to certify non-negativity.

Indeed, several non-SOS \germans{} for positive polynomials on compact sets following the general form~\eqref{eq:schmudgenform} have been proposed \citep{dressler2017positivstellensatz, dickinson2015extension, kuryatnikova2019copositive, chandrasekaran2016relative}. These are  {\em non-SOS Schm\"{u}dgen-type} \germans{} in which the polynomials $\sigma_\alpha$ in~\eqref{eq:schmudgenform}, instead of being SOS polynomials, belong to a different class of polynomials.
In particular, after adding appropriate redundant constraints to the definition of the set $S$, non-SOS Schm\"{u}dgen-type Positivstellens\"atze have been derived based on the classes of {\em sums of non-negative circuit} (SONC) polynomials~\citep{dressler2017positivstellensatz} and
 {\em sums of arithmetic-geometric mean exponential} (SAGE) functions~\citep[][Thm. 4.2]{chandrasekaran2016relative}.
In more generality, in~\citep[][Cor.~5]{kuryatnikova2019copositive}, it is shown that after adding appropriate redundant constraints to the definition of the set $S$, non-SOS Schm\"{u}dgen-type \germans{} can be derived based on any class of polynomials containing all positive constants. This condition is satisfied by SONC, \SAGElongtext{} (\SAGEtext{})~\citep{karaca2017repop, chandrasekaran2016relative},
{\em diagonally dominant SOS} (DSOS), and {\em scaled diagonally dominant SOS} (SDSOS) polynomials. Thus,
non-SOS Schm\"{u}dgen-type Positivstellens\"atze exist for all these classes of polynomials, and some of their numerical characteristics have been investigated~\citep[][]{kuryatnikova2019copositive, kuang2019alternative, dressler2019approach, dressler2020global}.

Schm\"{u}dgen's Positivstellensatz~\eqref{eq:schmudgenform} uses an exponential number of terms. For this reason Putinar's Positivstellensatz~\cite{putinar1993positive} is preferred in practice.
Under the assumption that the {\em quadratic module} associated with the polynomials $g_1,\dots,g_m$  is {\em Archimedean}~\citep[see, e.g.,][]{scheiderer2009positivity}, Putinar's Positivstellensatz~\cite{putinar1993positive} states that every polynomial $p$ that is positive on $S$, belongs to the quadratic module associated with $g_1,\dots,g_m$. That is, $p$ can be written in the form
\begin{equation}\label{eq:putinarform}
p(x) = \sigma_0(x) + \sum_{j=1}^m \sigma_j(x) g_j(x),
\end{equation}
where again, $\sigma_j$, for $j=0,1,\dots,m$, are SOS polynomials. Analogous to expression~\eqref{eq:schmudgenform}, expression~\eqref{eq:putinarform} makes the non-negativity of $p$ on $S$ evident. However, Putinar's Positivstellensatz~\eqref{eq:putinarform} uses a linear number of terms.

It is natural to ask if {\em non-SOS Putinar-type} Positivstellens\"atze do exists. That is, can Positivstellens\"atze following the general form~\eqref{eq:putinarform}, in which the polynomials~$\sigma_\alpha$ in~\eqref{eq:putinarform}, instead of being SOS polynomials, belong to a different class of  polynomials be derived? Indeed, when one specifically considers using SONC polynomials instead of SOS polynomials, the existence of the corresponding {\em SONC Putinar-type} Positivstellensatz has been posed as an open question in~\citep[][Sec. 6]{dressler2017positivstellensatz}.

Further motivation towards obtaining non-SOS Putinar-type Positivstellens\"atze can be gathered from the fact that simply replacing the SOS polynomials in~\eqref{eq:putinarform} by DSOS or SDSOS polynomials,
has been considered as a way to generate less computationally expensive approximation hierarchies to numerically approximate the solution of polynomial optimization problems~\citep[see, e.g.,][]{ahmadi2019dsos, kuang2016alternative, zheng2019sparse}. However, there is a lack of theoretical support
to guarantee the convergence of these approximation hierarchies; or whether a similar type of approximation hierarchy using the DSOS and SDSOS classes of polynomials or other classes of non-SOS polynomials can be guaranteed to converge. Thus far, the results obtained in this direction are negative~\citep{josz2017counterexample, ahmadi2017response, dressler2018optimization}.

In particular, it is known that replacing the SOS polynomials in~\eqref{eq:putinarform} by DSOS or SDSOS polynomials~\citep[as in][]{ahmadi2019dsos, kuang2016alternative, zheng2019sparse}, fails to lead to a {\em DSOS} or {\em SDSOS Putinar-type} Positivstellensatz~\citep[][]{josz2017counterexample, ahmadi2017response} (see Example~\ref{ex:counterSDSOS}). A key element of the proof of existing non-SOS Schm\"{u}dgen-type Positivstellensatz is the addition of redundant upper and lower bounds on the variables to the description of the compact set $S$~\citep{dressler2017positivstellensatz, kuryatnikova2019copositive, chandrasekaran2016relative}. Namely, the Positivstellens\"atze depend not only on the original polynomials $g_j$, $j=1,\dots,m$ describing $S$, but also on terms of the form $x_i-L_i$, $U_i -x_i$, where $L_i \le x_i \le U_i$, for $i = 1,\dots,n$.
When looking to obtain non-SOS Putinar-type \germans{}, it is reasonable
to add redundant lower and upper bounds on the variables to the set~$S$, as conjectured in~\citep{dressler2017positivstellensatz} for the SONC class of polynomials. However, it has been shown in~\citep[][Sec.~5]{dressler2018optimization}
 that such conjecture is false (see Example~\ref{ex:counterSONC}).

The aforementioned examples show that while the consideration of non-SOS classes of polynomials to certify non-negativity provides practical alternatives for the construction of non-SOS Putinar-type Positivstellens\"atze, it is a challenging problem to derive them. This is because
classes of polynomials such as  SONC, SDSOS, DSOS, and \SAGEtext{} polynomials, lack key properties of the class of SOS polynomials, such as being  closed under multiplication, which are key to prove the existence of certificates such as Putinar's Positivstellensatz via classical algebraic geometry tools~\citep[see, e.g.,][p. 544]{dressler2017positivstellensatz}. This shows that a novel approach, different to the ones followed in~\citep{ahmadi2019dsos, dressler2018optimization, dressler2017positivstellensatz} is needed to derive the desired non-SOS Putinar-type \germans{}.

 Notice that the constraints added to the compact set $S$ to obtain non-SOS \schmu-type \germans{} in~\citep{dressler2017positivstellensatz, kuryatnikova2019copositive, chandrasekaran2016relative} can be interpreted as
follows: On one hand, adding redundant inequalities to the description of $S$ maintains the set unchanged. On the other hand, these redundant inequalities enlarge the set of inequalities defining $S$, and therefore increase their expressive algebraic power, such that non-SOS expressions of the form~\eqref{eq:schmudgenform} capture all positive polynomials on $S$.
Thus, a natural question is whether one can add appropriate redundant inequalities to increase the expressive power of the inequalities defining $S$, such that  non-SOS expressions of the form~\eqref{eq:putinarform} capture all positive polynomials.

\subsection{Results Summary}

In positively answering the question above, the approach we apply is motivated by an interpretation of Putinar's Positivstellensatz as a Positivstellensatz for polynomials over the ball (see Theorem~\ref{thm:putinar}).
This interpretation allows us to circumvent the Archimedean requirements of Putinar's Positivstellensatz, and to apply the general ideas of~\citep{kuryatnikova2019copositive, pena2008exploiting} where Positivstellens\"atze over general semialgebraic sets are reduced to  Positivstellens\"atze over ``simpler'' sets.

Namely, let polynomials $g_1,\dots,g_m$ and $r > 0$ such that $S:=\{x \in \R^n: g_j(x) \ge 0, j=1,\dots,m\} \subseteq \{x\in\R^n:\|x\| \le r\}$ be given. Under mild assumptions on the class of non-negative polynomials $\cK$, we obtain non-SOS Putinar-type \germans{} (see Theorem~\ref{thm:semiSparse}); that is, if the polynomial $p$ is positive on $S$, then  there exist $2n+1$-variate polynomials $\alpha_j(y,z,u) \in \cK$ for $j=0,1$, and univariate polynomials $\rho_j(u) \in \cK$, $j=1,\dots,m$, such that
\begin{multline}\label{eq:PutFormIntro}
  p(x) = \alpha_0(x_1+r,\dots,x_n+r,r-x_1,\dots,r -x_n,r^2-\|x\|^2) \\
+ \alpha_1(x_1+r,\dots,x_n+r,r-x_1,\dots,r -x_n,r^2-\|x\|^2)(r^2-\|x\|^2)
+ \sum_{j=1}^m\rho_j(U_j - g_j(x))g_j(x),
\end{multline}
where $U_j$ is a large enough constant (see details next), for $j=1,\dots,m$.
To show that~\eqref{eq:PutFormIntro} makes the non-negativity of $p$ on $S$ evident, let us look in detail at  the right hand side of~\eqref{eq:PutFormIntro}. Let~$A(x)$  be the sum of the first two terms, and
$B(x) = \sum_{j=1}^m\rho_j(U_j - g_j(x))g_j(x)$. From $\alpha_0, \alpha_1 \in \cK$ (i.e., in particular being non-negative polynomials) it follows that $A(x)$ is non-negative on the ball $\{x\in\R^n:\|x\| \le r\}$. Further, from $\rho_j \in \cK$, and $g_j(x) \ge 0$ for all $x \in S$, for $j = 1,\dots,m$, it follows that $B(x) \ge 0$ on $S$. Since $S\subseteq \{x\in\R^n:\|x\| \le r\}$, then $p(x) = A(x) + B(x)$ is non-negative on $S$.

This partition of the right hand side~\eqref{eq:PutFormIntro} also allows us to explain the high level idea of our prof of the existence of the \germans{}~\eqref{eq:PutFormIntro}. The proof consist of two steps. The first one is to show that there exists~$\rho_j \in \cK$ for $j=1,\dots,m$, such that $p(x) - B(x) > 0$ on $\{x\in\R^n:\|x\| \le r\}$.
We will show (see Lemma~\ref{lem:StoB}), that this is the case when the class $\cK$ contains the class of sums of monomials squares (SOMS) polynomials (cf., Definition~\ref{def:soms}). The second step is to show that every positive polynomial on $\{x\in \R^n:\|x\| \le r\}$ can be written in of the form of $A(x)$. This will imply that $p(x) - B(x) = A(x)$, which implies~\eqref{eq:PutFormIntro}. We will show (see Lemma~\ref{lem:KputBall}) that it is always possible to obtain a Positivstellensatz on the ball of the form $A(x)$ where~$\alpha_0$, $\alpha_1$ are SOMS polynomials. Thus, the only assumption  on the class of polynomials~$\cK$ required for the \germans{}~\eqref{eq:PutFormIntro} to exist is that it contains contains all SOMS polynomials (see Theorem~\ref{thm:semiSparse}).

Note that as long as $U_j \ge \max\{g_j(x): \|x\| \le r\}$, the arguments above, about the right hand side of~\eqref{eq:PutFormIntro}, hold even for some not necessarily non-negative classes of polynomials, such as \SAGEtext{} polynomials that are non-negative only in the positive orthant. This follows from the fact that in this case, the polynomials $\alpha_0, \alpha_1$, and $\rho_j$ for $j=1,\dots,m$, are evaluated on entries that are always non-negative when $x \in S$. Thus, as long as $U_j$ is large enough for $j=1,\dots,m$ (see Proposition~\ref{prop:boundU} for an explicit bound), the condition of $\cK$ being a class of non-negative polynomials is not needed (see Theorem~\ref{thm:semiSparse} and the discussion that follows). This fact will allow us to capture the case in which \SAGEtext{} polynomials are used as the class $\cK$ in the proposed \germans{} (see Remark~\ref{rem:RHS2}).

The classes of SONC, SDSOS, DSOS and \SAGEtext{} polynomials contain all SOMS polynomials (see Proposition~\ref{prop:relationKs}). Thus, any of these classes of non-SOS polynomials can be used as the class of polynomials $\cK$ in which the \germans{}~\eqref{eq:PutFormIntro} is based (see Corollary~\ref{cor:semiSparse}). In particular, we provide the SONC Putinar-type \german{} seeked in~\citep[][Sec.~6]{dressler2017positivstellensatz}.

Another key property of the \germans{}~\eqref{eq:PutFormIntro} is their inherent sparsity characteristics, due to the use of the univariate polynomials $\rho_j$. Such sparsity structure can be further enhanced when the polynomial $p$,  whose non-negativity is being certified, as well as the polynomials defining the underlying semialgebraic set $S$, are sparse.
Deriving \germans{} that take advantage of both the sparsity of $p$ and the polynomials defining $S$ has been the focus of a wealth of research work. For example, consider the earlier work in~\citep{parrilo2002explicit, pena2008exploiting, laurent2007semidefinite,waki2006sums}, and the more recent work in~\citep{lasserre2006convergent, mai2020sparse,wang2020cs,wang2021tssos}. This latter work focuses on exploiting {\em correlative and term sparsity} to derive term and correlative sparse versions of SOS \germans{} such as Putinar's Positivstellensatz (in~\citep{lasserre2006convergent,wang2020cs,wang2021tssos}); and Reznick's and Putinar-Vasilescu's~\citep{putinar1999solving} Positivstellensatz (in~\citep{mai2020sparse}). The \germans~\eqref{eq:PutFormIntro} guarantees the existence of non-SOS Putinar-type \germans{} which are term sparse, without any assumption on the polynomial $p$ or the polynomials defining $S$ (see Theorem~\ref{thm:semiSparse} and Corollary~\ref{cor:semiSparse}). Moreover (see Theorem~\ref{thm:mainSparse}), we
show how to take advantage of correlative sparsity to derive sparse versions of the non-SOS Putinar-type Positivstellensatz~\eqref{eq:PutFormIntro}. Note that in contrast with the mentioned previous related work, our results show how to exploit sparsity in a more general and novel setting in which SOS polynomials are not necessarily used as the base class of polynomials to certify non-negativity, but instead other classes of polynomials such as SONC, SDSOS, DSOS, and \SAGEtext{} polynomials could be chosen for this purpose (see Corollary~\ref{cor:mainSparse}).

Although the focus throughout the article is the use of non-SOS classes of polynomials to derive new Putinar-type \germans{}; it is clear that the  class of SOS polynomials contains all SOMS polynomials. That is, the class of SOS polynomials can be used as the base class in the Putinar-type \germans{}~\eqref{eq:PutFormIntro}.
Thus, as a byproduct, we obtain an alternative way to derive versions of known sparse SOS Putinar-type \germans{}, with additional sparsity characteristics, from the proposed non-SOS Putinar-type \germans{} (see Section~\ref{sec:SOS}).

\subsection{Organization}

The rest of the paper is structured as follows. In Section~\ref{sec:notation}, we introduce the main notation and concepts used throughout the article, including well-known SOS \germans{} and classes of polynomials that have been used to certify non-negativity. In Section~\ref{sec:about}, we discuss in more detail recent results regarding non-SOS Schm\"udgen-type \germans{}, as well as the motivation behind and challenges of obtaining non-SOS Putinar-type \germans{}. In Section~\ref{sec:NewputinarType}, we present the main results of the article. Namely,
in Section~\ref{sec:semisparse}, we derive Putinar-type \germans{} for a wide class of non-SOS classes of polynomials (see Theorem~\ref{thm:semiSparse} and Corollary~\ref{cor:semiSparse}), including  SONC, SDSOS, DSOS, and \SAGEtext{} polynomials. These \germans{} have inherent sparsity characteristics that are further exploited in Section~\ref{sec:fullysparse}, when the sparsity structure of both the polynomial whose non-negativity is being certified and the polynomials defining the semialgebraic set of interest are known (see Theorem~\ref{thm:mainSparse} and Corollary~\ref{cor:mainSparse}). In Section~\ref{sec:examples} we present some examples that illustrate how the derived \germans{} resolve the problem of certifying the non-negativity of some particular polynomials considered in the literature, using non-SOS Putinar-type \germans{}. Further, we provide some examples that
 illustrate the effect that the use of different classes of polynomials and sparsity information has in the computational effort needed to compute the proposed \germans{}. For the purpose of clarity, the proofs of a number of the results are delayed to Section~\ref{sec:proofs}. In particular, we provide the proof of Theorem~\ref{thm:semiSparse} in Section~\ref{sec:proofsemisparse} and the proof of Theorem~\ref{thm:mainSparse} in Section~\ref{sec:prooffullysparse}. Also, in Section~\ref{sec:SOS}, we consider the implications of our results for the particular case in which the class of SOS polynomials is used to certify non-negativity.
In Section~\ref{sec:end}, we finish with some conclusions and directions for future work.

\section{Preliminaries}
\label{sec:notation}

For any $r > 0$, let $\B_{r}:= \{x \in \R^n: \|x\| \le r\}$ be the ball of radius $r$ centered at the origin.
Denote by $\R[x]:=\R[x_1,\dots,x_n]$ the set of $n$-variate polynomials with real coefficients. Also, let $e$ denote the all-ones vector of appropriate dimensions.

We will use the following notation to ease the presentation. Given polynomials $g_1,\dots,g_m \in \R[x]$ and~$\alpha \in \N^{m}$, let $g := [g_1,\dots,g_m]\tr$ be the vector of polynomials whose entries are $g_1,\dots,g_m$  and
$g^{\alpha} := \prod_{j=1}^m g_j^{\alpha_j}$. In particular, $x^\alpha = \prod_{j=1}^m x_j^{\alpha_j}$.

A polynomial $p \in \R[x]$ is called non-negative (resp. positive) on $S \subseteq \R^n$ if $p(x)\geq 0$ (resp. $p(x)>0$) for all $x\in S$. For the special case in which $p \in \R[x]$ is non-negative (resp. positive) on $S = \R^n$, we will simply say that the polynomial $p$ is non-negative (resp. positive). Further, we say that $\cK \subset \R[x]$ is a class of non-negative polynomials if every $p \in \cK$ is non-negative. For instance, the class of {\em sum of squares} (SOS) polynomials is a class of non-negative polynomials.

\begin{definition}[SOS polynomials]
\label{def:SOS}
A polynomial $p\in \R[x]$ is an SOS if $p=\sum_{l=1}^sq_l^2$ for some $s \in \N$, where $q_l\in \R[x]$ for $l=1,\dots,s$. We will denote by $\KSOS$ the class of SOS polynomials; that is, $\KSOS = \{ p \in \R[x]: p \text{ is an SOS polynomial}\}$.
\end{definition}

SOS polynomials can be written in {\em Gram matrix} form where the Gram matrix of coefficients is a positive semidefinite matrix. As a result, membership in $\KSOS$ can be tested by solving a semidefinite programming (SDP) problem~\citep[see, e.g.,][]{BlekPT13}.

Given $g_1,\dots,g_m\in \R[x]$, and a class of polynomials $\cK \subseteq \R[x]$, let us also define
\[
\cP_{\cK}(g)= \left \{\sum_{\alpha \in \N^m_d}s_{\alpha}g^{\alpha}:s_{\alpha} \in \cK,\, d=0,1,\dots \right \}
\]
and
\[
\cM_{\cK}(g)=\left \{s_0 + \sum_{j=1}^m s_jg_j: s_j \in \cK, j=0,\dots,m \right \}.
\]
Note that for any class of non-negative polynomials $\cK$, $p \in  \cP_{\cK}(g) \supseteq \cM_{\cK}(g)$, implies that $p$ is a non-negative polynomial on the semialgebraic set
\[
S_g:=\{x\in \R^n:g_j(x)\ge 0, j=1,\dots,m\}.
\]
When $\cK = \KSOS$, the sets $\cP_{\KSOS}(g)$ and $\cM_{\KSOS}(g)$ are respectively known~\citep[see, e.g.,][]{marshall2003approximating} as the \emph{preorder} and the \emph{quadratic module}  associated with $g$. Also, when $\cK = \R_+$; that is, $\cK$ is the class of non-negative constant polynomials, $\cP_{\R_+}(g)$ is the preprime associated with the semialgebraic set~$S_g$, and for any class of non-negative polynomials $\cK \in \R[x]$, $\cP_{\cK}(g)$ is a module over $\cP_{\R_+}(g)$~\cite[see, e.g.,][]{marshall2003approximating}.

Schm\"{u}dgen's Positivstellensatz~\cite{schmudgen1991k} states that all positive polynomials on a semialgebraic compact set belong to the preorder associated with the polynomials defining the set.

\begin{theorem}[Schm\"{u}dgen's Positivstellensatz~\cite{schmudgen1991k}]
\label{thm:schmudgen}
Let polynomials $g_1,\dots,g_m \in \R[x]$ be given such that $S_g = \{x \in \R^n: g_j(x) \ge 0, j=1,\dots,m\}$ is compact. Then $p > 0$ on $S_g$
implies $p \in \cP_{\KSOS}(g)$.
\end{theorem}

Schm\"{u}dgen's Positivstellensatz (Theorem~\ref{thm:schmudgen}) is constructed using an exponential number, $2^{m}$, of $n$-variate $\SOS$ polynomials (terms). From a computational point of view, these certificates are therefore expensive~\citep[see, e.g.,][]{lasserre2006convergent}. When  $\cM_{\KSOS}(g)$ is {\em Archimedean}, Putinar's Positivstellensatz~\cite{putinar1993positive} shows that every polynomial positive on the (compact) semialgebraic set $S_g$ belongs to the quadratic module $\cM_{\KSOS}(g)$; and therefore, it offers a more practical alternative to Schm\"{u}dgen's Positivstellensatz as it only uses a linear number, $m +1$, of $n$-variate $\SOS$ polynomials (terms). Putinar's Positivstellensatz is one of the most important tools in the area of polynomial optimization. In fact, it is the building block of Lasserre's seminal work in polynomial optimization~\citep[see, e.g.,][]{anjo12, BlekPT13, lasserre2009moments}. Theorem~\ref{thm:putinar} below is an equivalent form of Putinar's Positivstellensatz~\cite{putinar1993positive}.

\begin{theorem}[Putinar's Positivstellensatz]
\label{thm:putinar}
 Let $r > 0$ and polynomials $g_1,\dots,g_m \in \R[x]$ be given. Then $p(x) > 0$ on $\{x \in \B_{r}: g_j(x) \ge 0, j=1,\dots,m\}$
implies $p(x)\in \cM_{\KSOS}(g,r^2 - \|x\|^2)$.
\end{theorem}

Note that since  $\cM_{\KSOS}(g,r^2 - \|x\|^2)$ is  Archimedean, Theorem~\ref{thm:putinar} is implied by Putinar's Positivstellensatz~\cite{putinar1993positive}. Further, if $\cM_{\KSOS}(g)$ is Archimedean, taking $r > 0$ such that $r^2 - \|x\|^2 \in \cM_{\KSOS}(g)$ we obtain $\cM_{\KSOS}(g(x), r^2 - \|x\|^2)  \subseteq \cM_{\KSOS}(g)$ and therefore Theorem~\ref{thm:putinar} implies the original Putinar's Positivstellensatz for the set $S_g$.

Theorem~\ref{thm:putinar} allows us to circumvent the need for the Archimedean assumption by working with the set $\{x \in \B_{r}: g(x) \ge 0\}$. In practical terms we are trading the Archimedean assumption, which is equivalent to having a Putinar \german{} for the positivity of  $r^2 - \|x\|^2$ in $S_g$,  with the (weaker) assumption of knowing that $r^2 - \|x\|^2$ is non-negative on $S_g$.

A generalization of Theorem~\ref{thm:putinar} (see Lemma~\ref{lem:StoB}) will be one of the main building blocks to show the existence of \germans{} of the general form of Putinar's \german{} in which, loosely speaking, the class of SOS polynomials $\KSOS$ is replaced by any class of polynomials containing all the sums of monomial squares  (see Theorem~\ref{thm:semiSparse}), and in particular, by some well known classes polynomials (see Corollary~\ref{thm:semiSparse}).

\begin{myremark}
Analogous to the case of SOS polynomials (see, Definition~\ref{def:SOS}), throughout the article, given a type of polynomials $(\cdot)$, we will use $\cK_{(\cdot)}$ to indicate the class of polynomials of that type.
\end{myremark}

Next, we formally introduce these latter classes of polynomials and discuss some of their properties and relationships. In line with the theme of the results in this article, these classes of polynomials have been recently considered or introduced in the areas of algebraic geometry and polynomial optimization, as an alternative to the use of SOS polynomials.

We begin by introducing the simple class of {\em sum of monomial squares} (SOMS) polynomials that plays a central role in our results.

\begin{definition}[SOMS polynomials]
\label{def:soms}
A polynomial $p \in \R[x]$ is a {\em SOMS} if $p(x) = \sum_{\alpha \in A} \lambda_{\alpha} x^{2\alpha}$ for some $A \subseteq \N^n$ and $\lambda_{\alpha}  \ge 0$ for all $\alpha \in A$.
\end{definition}
Clearly, $\KSOMS \subset \KSOS$; in particular, SOMS polynomials are non-negative. Also, membership in $\KSOMS$ can be tested by solving a linear programming (LP) problem.

Next, consider the (non-SOS) classes of {\em diagonally dominant SOS} (DSOS) and {\em scaled diagonally dominant SOS} (SDSOS) polynomials~\citep{ahmadi2019dsos, fidalgo2011positive}.

\begin{definition}[DSOS and SDSOS polynomials]
A polynomial $p \in \R[x]$ is a SDSOS polynomial if
$
p(x)=\sum_{l=1}^k (a_{l}x^{\alpha_{l}}+b_{l}x^{\beta_{l}})^2,
$
for some $k \in \N$, and some $a_{l}, b_{l} \in \R$ and $\alpha_{l}, \beta_{l} \in \N^n$, for $l=1,\dots,k$. If further $a_{l} = b_{l}$ or $a_{l} = -b_{l}$ for each $l=1,\dots,s$, then $p$ is a $\DSOS$ polynomial.
\end{definition}

These polynomials are also referred as sums of binomial squares~\citep[see, e.g.,][]{wang2018nonnegative}.
Note that $\KDSOS \subset \KSDSOS \subset \KSOS$. Also, as their name indicates, SDSOS (resp. DSOS) polynomials can be written in Gram matrix form where the Gram matrix of coefficients is a scaled diagonally dominant (resp. diagonally dominant) matrix. As a result, membership in $\KSDSOS$ (resp. $\KDSOS$) can be tested by solving a second-order cone programming (SOCP) problem (resp. an LP problem)~\citep{ahmadi2019dsos}. In particular, in some polynomial optimization applications it is computationally more advantageous to test for membership in $\KDSOS$ than in $\KSOS$~\citep[see, e.g.,][]{ahmadi2019dsos, kuang2016alternative, zheng2019sparse}.

Now consider the class of {\em circuit polynomials}~\citep{pantea2012global,  iliman2016amoebas}.
A polynomial $p \in \R[x]$ is a {\em circuit polynomial} if
$p(x)=\sum_{\alpha \in A}c_{\alpha}x^{\alpha} - dx^{\beta}$,
where $c_{\alpha} > 0$,  for all $\alpha \in A$, $A \subseteq (2\N)^{n}$ is the vertex set of a simplex, and~$\beta$ lies in the interior of this simplex.
Circuit polynomials that are non-negative were characterized in~\citep[][Thm. 3.8]{iliman2016amoebas} based on their {\em circuit number}. This characterization leads to the definition of  {\em sum of non-negative circuit} (SONC) polynomials considered in~\citep[see, e.g.,][]{wang2018nonnegative,dressler2017positivstellensatz, dressler2018optimization, iliman2016amoebas}, whose membership can be tested using geometric programming~\citep{dressler2019approach, iliman2016lower, papp2019duality}, or SOCP~\citep{magron2020sonc, wang2020second}.

\begin{definition}[SONC polynomials]
A polynomial $p \in \R[x]$ is a SONC polynomial if
$p(x)=\sum_{l=1}^s q_l(x)$ for some $s \in \N$, where $q_l$ is a non-negative circuit polynomial for $l=1,\dots,s$.
\end{definition}

The last class of polynomials we introduce, is the class of {\em \SAGElongtext{}} (\SAGEtext{}) polynomials~\citep[see, e.g.,][]{karaca2017repop, chandrasekaran2016relative}, which are closely related to the SONC polynomials introduced above~\citep[see, e.g.,][for details]{wang2018nonnegative,dressler2017positivstellensatz, murray2021applications}.

\begin{definition}[\SAGEtext{} polynomials]
A polynomial $p\in \R[x]$ is a \SAGEtext{} polynomial if $p(x)=\sum_{l=1}^s q_l(x)$ for some $s \in \N$, where $q_l$ is an \AGElongtext{} (\AGEtext{}) polynomial for $l=1,\dots,s$. In turn, $q\in \R[x]$ is an \AGEtext{} polynomial if $q$ is non-negative on the non-negative orthant (i.e., $q(x) \ge 0$ for all $x \in \R^n_+$) and has at most one negative coefficient.
\end{definition}
Membership in class of \SAGEtext{} polynomials can be tested using relative entropy programming (REP)~\citep{karaca2017repop, chandrasekaran2016relative, murray2018newton} and geometric programming~\citep{ghasemi2012lower}.

Next, we state a relationship between the classes of polynomials discussed above that will be key for the results presented later.

\begin{proposition}\label{prop:relationKs}
We have the following relations
\begin{enumerate}[label = (\roman*)]
\item
$\KSOMS \subseteq \KDSOS \subseteq \KSDSOS$,
\label{it:relation1}
  \item  $\KSDSOS = \KSOS \cap \KSONC$,
 \label{it:relation2}
 \item $\KSONC \subseteq \KSAGE$.
 \label{it:relation3}
\end{enumerate}
\end{proposition}

\begin{proof}
Statement~\ref{it:relation1} clearly follows from the definitions of SOMS, DSOS, and SDSOS polynomials.  The inclusion  $\KSOS \cap \KSONC \subseteq \KSDSOS$ is shown in~\citep[][Thm. 5.2]{iliman2016amoebas}. The inclusion $\KSDSOS \subseteq \KSOS \cap \KSONC$ follows from the inclusions $\KSDSOS \subseteq \KSONC$ shown in~\citep[see, e.g.,][Lem. 4.1]{kurpisz2019new}; and $\KSDSOS \subseteq \KSOS$, which follows from the definition of SDSOS and SOS polynomials. This proves statement~\ref{it:relation2}.
Statement~\ref{it:relation3} follows from the definitions of SONC and \SAGEtext{} polynomials; that is, any non-negative circuit polynomial is an \AGEtext{} polynomial.
\end{proof}

We finish this section by presenting a bound for the value of a polynomial over a ball that will be relevant for the results presented later. A polynomial $g \in \R[x]$ of degree $d$ can be written as $g(x)=\sum_{\{\alpha\in\mathbb{N}^n: e\tr \alpha\leq d\}}g_{\alpha}x^{\alpha}$.
Let ${\|g\|}_2^2:=\sum_{\{\alpha\in\mathbb{N}^n: e\tr \alpha\leq d\}} \frac{\alpha_1!\cdots\alpha_n!(d-e\tr \alpha)!}{d!}g_{\alpha}^2$.

\begin{proposition}\label{prop:boundU}
Let $g \in \R[x]$ be a polynomial of degree $d \ge 1$, and $r>0$. Then,
\[
\max\{g(x): x \in \B_r\} \le {\|g\|}_2 (1+r^2)^{d/2} .
\]
\end{proposition}

\begin{proof}
Let $x_0:= 1$.
For each $i_1,i_2,\dots,i_d \in \{0,\dots,n\}$, let $\alpha = (\alpha_1, \dots, \alpha_n)$ be such that
$x_{i_1}x_{i_2} \cdots x_{i_d} = x_1^{\alpha_1} \cdots x_n^{\alpha_n}$, and let $c_{i_1,\dots,i_d} = \frac{\alpha_1!\cdots\alpha_n!(d-e\tr \alpha)!}{d!}g_\alpha$. Then,
we have that \[g(x) = \sum_{i_1 =0}^n\sum_{i_2 =0}^n \cdots \sum_{i_d = 0}^n  c_{i_1,\dots,i_d} x_{i_1}x_{i_2}...x_{i_d}.\]
Using Cauchy-Schwarz, it follows that
\begin{align*}
g(x)^2  & \le \sum_{i_1 =0}^n \sum_{i_2 =0}^n\cdots\sum_{i_d = 0}^n  c_{i_1,\dots,i_d}^2    \sum_{i_1 =0}^n \sum_{i_2 =0}^n\cdots\sum_{i_d = 0}^n   (x_{i_1}x_{i_2}\cdots x_{i_d})^2 \\
            &  = \sum_{i_1 =0}^n \sum_{i_2 =0}^n \cdots\sum_{i_d = 0}^n  c_{i_1,\dots,i_d}^2     \sum_{i_1 =0}^n x_{i_1}^2   \sum_{i_2 =0}^n  x_{i_2}^2  \cdots \sum_{i_d = 0}^n  x_{i_d}^2\\
             & = \sum_{i_1 =0}^n \sum_{i_2 =0}^n \cdots\sum_{i_d = 0}^n  c_{i_1,\dots,i_d}^2   (1+\|x\|^2)^d.
\end{align*}
Also notice that
\[
\sum_{i_1 =0}^n \sum_{i_2 =0}^n \cdots\sum_{i_d = 0}^n  c_{i_1,\dots,i_d}^2  =   \sum_{\{\alpha\in\mathbb{N}^n: e\tr \alpha\leq d\}} \frac{\alpha_1!\cdots\alpha_n!(d-e\tr \alpha)!}{d!}g_{\alpha}^2 = {\|g\|}_2^2.
\]
Therefore, if  $\|x\| \le r$, we obtain  that $g(x) \le {\|g\|}_2 (1+r^2)^{d/2}$.  
\end{proof}

\section{About non-SOS Positivstellens\"atze}
\label{sec:about}

As mentioned earlier, it is computationally expensive to check membership in $\KSOS$ to construct \germans{} such as Schm\"{u}dgen's (Theorem~\ref{thm:schmudgen}) and Putinar's (Theorem~\ref{thm:putinar}) \german{}. This fact has led to the proposal of alternative classes of polynomials $\cK \neq \KSOS$; such $\cK = \KDSOS$, $\cK =\KSDSOS$, $\cK =\KSONC$, and $\cK = \KSAGE$, to be used as the base class of polynomials to construct \germans{} for polynomials over compact sets; that is, {\em non-SOS Positivstellens\"atze}. In particular,
{\em non-SOS Schm\"{u}dgen-type Positivstellensatz} have been recently derived in~\citep{dressler2017positivstellensatz, dickinson2015extension, kuryatnikova2019copositive, chandrasekaran2016relative}; which follow the general form~\eqref{eq:schmudgenform}. Namely, for some classes of polynomials $\cK \neq \KSOS$ and~$g \in \R[x]^m$ such that $S_g$ is compact, these results show that every positive polynomial on $S_g$ belongs to $\cP_{\cK}(\bar g)$ where $\bar g$ is obtained from $g$ by adding appropriate redundant constraints such that $S_g = S_{\bar g}$. More specifically, $\bar g := (g,x-L,U-x)$. where $L,U \in \R^n$  are  such that $S_g \subseteq \{x \in \R^n:L \le x \le U\}$.

More specifically, non-SOS Schm\"{u}dgen-type \germans{} have been obtained by using the classes $\KSONC$ and $\KSAGE$ using an abstract \german{} from real algebraic geometry that requires appropriate Archimedean assumptions~\cite{marshall2008positive}. Assume $S_g$ is compact and $p \in \R[x]$ is positive on $S_g$, Following~\cite{marshall2008positive}, \citet[][]{dressler2017positivstellensatz}
show that $p\in \cP_{\KSONC}(\bar{g})$. Similarly,~\citet[][Thm. 4.2]{chandrasekaran2016relative},
show that $p\in \cP_{\KSAGE}(\bar{g})$. In both works the valid inequalities $x-L$, $U-x$, giving explicit upper and lower bounds, are added to $g$ to ensure the required Archimedean assumptions from~\cite{marshall2008positive}.
More recently, in~\citep[][Cor.~6]{kuryatnikova2019copositive}, it is shown that $p\in \cP_{\R_+}(\bar{g}) \subseteq \cP_{\KDSOS}(\bar{g}) \subseteq \cP_{\KSDSOS}(\bar{g})$, showing that these non-SOS Schm\"{u}dgen-type \germans{} can also be obtained using the classes of $\KDSOS$ and $\KSDSOS$ as the base class of polynomials to certify non-negativity.

All these non-SOS Schm\"{u}dgen-type Positivstellensatz follow from more general \germans{} recently introduced in~\citep[][Cor.~5]{kuryatnikova2019copositive}, which are obtained, similarly to~\cite{marshall2008positive}, by reducing the problem of certifying non-negativity over a compact semialgebraic set to the case when the semialgebraic set is the simplex, and then applying Polya's \german{}~\citep[see, e.g.,][]{hardy1952inequalities}. This result is presented below in a form that is amenable for use in what follows.

\begin{theorem}[$\cK$-Schm\"{u}dgen \germans{}~{\citep[][Cor.~5]{kuryatnikova2019copositive}}]\label{Thm:Kschmudgen} Let $\cK \subseteq \R[x]$ be a class of polynomials such that $\R_+ \subseteq \cK$. Let polynomials $g_1,\dots,g_m \in \R[x]$ be given such that
$S_g = \{x \in \R^n : g_j(x) \ge 0, j=1,\dots,m\}$ is compact. Let $L,U \in \R^n$ be such that $S_g \subseteq \{x \in \R^n:L \le x \le U\}$. Then,
$p>0$ on $S_g$ implies $p\in \cP_{\cK}(g,x-L,U-x)$.
\end{theorem}
\begin{proof}
Let $M= e\tr U$, $\bar{g} := [g, x-L, U-x]\tr$, and $\widehat{g} := [g(x), x-L, M-e\tr x]\tr$. Note that $S_g = S_{\bar{g}} = S_{\widehat{g}}$. Then the result follows from~\citep[][Cor.~5]{kuryatnikova2019copositive} after noticing that $M-e\tr x = \sum_{i=1}^n (U_i - x_i)$ implies that  $\cP_{\cK}(\widehat{g}) \subseteq \cP_{\cK}(\bar{g})$.
\end{proof}
Note that since $\R_+ \subset \KDSOS \subset \KSDSOS\subset \KSONC$, and $\R^+  \subseteq \KSAGE$
the non-SOS Schmudgen-type Positivstellensatz for SONC, SDSOS, DSOS,  and \SAGEtext{} polynomials
follows from the $\cK$-Schm\"{u}dgen Positivstellensatz (Theorem~\ref{Thm:Kschmudgen}).

In light of these results, it is natural to ask whether a $\cK$-Putinar Positivstellensatz  for non-SOS classes $\cK$ of polynomials, analogous to the $\cK$-Schm\"{u}dgen Positivstellensatz (Theorem~\ref{Thm:Kschmudgen}), can be proved to exist.
The authors in~\citep[][p. 544, and Sec. 6]{dressler2017positivstellensatz} had posed the question for the case of SONC polynomials. But, proving the existence of {\em non-SOS Putinar-type} \germans{}; which follow the general form of~\ref{eq:putinarform}, but use alternative classes of polynomials $\cK \neq \KSOS$ as the basis to certify non-negativity, has been shown to be a challenging task. One main reason for this, is the fact that
the class of SOS polynomials is closed under multiplication, while other proposed classes of non-SOS polynomias such as DSOS, SDSOS, SONC, and \SAGEtext{} are not~\citep[see, e.g.,][p. 544]{dressler2017positivstellensatz}.

For example, it is known that a simple extension of Putinar's Positivstellensatz fails for non-SOS classes $\cK$ of non-negative polynomials over compact semialgebraic sets, even when the associated quadratic module is Archimedean. More specifically, when the non-SOS class of polynomials in consideration is the class of DSOS or SDSOS polynomials, this failure is illustrated by Example~\ref{ex:counterSDSOS} (see \citep{josz2017counterexample} and~\citep[][Sec.~2.2]{ahmadi2017response}).

\begin{myexample}[Simple DSOS/SDSOS Putinar Positivstellensatz fails~{\citep{josz2017counterexample}}]
\label{ex:counterSDSOS}
Let $g(x_1,x_2) = 1-x_1^2-x_2^2$ and $S_g = \{ (x_1,x_2) \in \R^2: 1-x_1^2-x_2^2 \ge 0)\}$.
The quadratic module generated by $g(x_1,x_2)$ is Archimedean since $1-\|(x_1,x_2)\|^2 = g(x_1,x_2)$. Thus, any positive polynomial on $S_g$ belongs to $\cM_{\KSOS}(g)$. Let  $p(x_1,x_2)=(-2+x_1+x_2)^2$, for any $(x_1,x_2) \in S_g$, $x_1 + x_2\le \sqrt{2} $. Thus $p(x) \ge 6 - 4\sqrt 2 > 0$ for all $x \in S_g$. However, in~\citep[eq. (4.5)]{josz2017counterexample}, it is shown that there are no SDSOS (and by extension not DSOS) polynomials $\sigma_0(x_1,x_2)$, and $\sigma_1(x_1,x_2)$, such that
$p(x_1,x_2)=\sigma_0(x_1,x_2)+\sigma_1(x_1,x_2)(1-x_1^2-x_2^2)$. That is, $p \not \in \cM_{\KSDSOS}(g)$
\end{myexample}

In light of the $\cK$-Schm\"{u}dgen Positivstellensatz (Theorem~\ref{Thm:Kschmudgen}), it is natural to consider if the desired Putinar-type Positivstellensatz can be obtained using $\cM_{\cK}(g, x-L, U-x)$, for non-SOS classes $\cK$ of polynomials~\citep[see, e.g.,][Sec. 6]{dressler2017positivstellensatz}. However, as
illustrated by Example~\ref{ex:counterSONC}~\citep[][Sec.~5]{dressler2018optimization}, this type of proposed Putinar-type Positivstellensatz fails for SONC polynomials.

\begin{myexample}[Proposed SONC Putinar-type Positivstellensatz fails~{\citep[][Sec.~5]{dressler2018optimization}}]
\label{ex:counterSONC} Let $n > 2$. Let $g_i(x) = 1-x_i^2, g_{n+i}(x) = x_i^2-1$, $i=1,\dots,n$. The quadratic module generated by $g$ is Archimedean,  as $n-\|x\|^2 = \sum_{i=1}^n g_i(x)$. Thus, any positive polynomial on $S_g= \{\pm 1\}^n$, belongs to $\cM_{\KSOS}(g)$. Let  $a>\tfrac{2^n-1}{2^{n-2}-1}$, and let
$p(x)=1+(a-1)\smash{\prod_{i=1}^n}\frac{x_i+1}{2}$.
Now let $c$ be such that $c_i \ge 1$ for $i=1,\dots,n$. Set $\bar{g} = [g, c-x, x-c]$. It follows that $S_{\bar{g}} = S_g = \{\pm 1\}^n$. Notice that $p(e) = a >0$, and for any $x \in \{\pm 1\}^n \setminus \{e\}$, $p(x) = 1$. Thus, $p(x) > 0$ on $\{\pm 1\}^n$. However, in~\citep[][Thm. 5.1]{dressler2018optimization},  it is shown that
$p(x) \notin \cM_{\KSONC}(g,x-c,c-x)$.

\end{myexample}
\begin{myremark} Example~\ref{ex:counterSONC} can be extended to cover the case of more general upper and lower bounds. Given  $L, U \in \R^n$ satisfying $S_g \subseteq \{x \in \R^n: L \le x \le U\}$, and $r \in \R$ satisfying $S_g \subseteq \{x: \|x\| \le r\}$, let $c_i = \max\{|L_i|,U_i\}$. Then $c_i \ge 1$ and using
$x - L_i = x -c_i + (c_i - |L_i|)$, $U_i-x = c_i-x + (c_i-U_i)$, for $i=1,\dots,n$, and $r^2 - \|x\|^2  = \sum_{i=1}^{n} g_i(x) + (r^2 - n)$
we have  that $\cM_{\KSONC}(g,x-L,U-x,r^2 - \|x\|^2) \subseteq \cM_{\KSONC}(g,x-c,c-x)$ and thus $p \notin \cM_{\KSONC}(g,x-L,U-x,r^2 - \|x\|^2)$ either.
\end{myremark}

In the next section, we show the existence of non-SOS Putinar-type \germans{} for non-SOS classes of polynomials  $\cK$ satisfying the mild assumption $\KSOMS \subseteq \cK$ (Theorem~\ref{thm:semiSparse}). Our result applies to a wide range of classes of polynomials $\cK$, including $\KDSOS$, $\KSDSOS$, $\KSONC$, and $\KSAGE$ (see Corollary~\ref{cor:semiSparse}).

As $\KSOS$ also satisfies these mild assumptions (see Section~\ref{sec:SOS}), it is easy to see that Theorem~\ref{thm:semiSparse} implies Theorem~\ref{thm:putinar}. Indeed it is enough to notice that $\KSOS$ is closed under composition with any polynomial (see Lemma~\ref{lem:simplifySOS} and Section~\ref{sec:SOS}). Notice that because in general, classes of polynomials do not have the property of being closed under such composition, those compositions appear explicitly in~\eqref{eq:PutFormIntro}. Our results (and methods) show that $\cK$ being closed under products is not an essential property to obtain non-SOS Putinar-type \germans{}.

\section{Sparse non-SOS Putinar-type Positivstellens\"atze}
\label{sec:NewputinarType}

The only requirement on the class~$\cK$ of polynomials  in the $\cK$-Schm\"{u}dgen's \germans{} (Theo\-rem~\ref{Thm:Kschmudgen}) is for $\cK$ to contain all constant non-negative polynomials; that is,  $\R_+ \subseteq \cK$. Obtaining an analogous non-SOS Putinar-type Positivstellensatz (Theorem~\ref{thm:semiSparse}, next) requires a bit more. Namely, $\cK$ must contain $\KSOMS$ (see Definition~\ref{def:soms}). As stated in Proposition~\ref{prop:relationKs}, this condition is satisfied by the classes of DSOS, SDSOS, SONC, and \SAGEtext{} polynomials.

\begin{myremark}
Results through the Section require an upper bound on the expression $\max\{q(x): x \in \B_r\}$ for some~$r >0$, and~$q \in \R[x]$. Proposition~\ref{prop:boundU} provides a constructive upper bound for this maximum.
 \end{myremark}

\subsection{Semi-sparse non-SOS Putinar-type Positivstellens\"atze}
\label{sec:semisparse}

We are ready to state the first main result of the article. Namely, non-SOS Putinar-type \germans{} for general classes of polynomials satisfying a couple of mild conditions, that exhibits inherent sparsity characteristics, and whose proof is delayed to Section~\ref{sec:proofsemisparse} for ease of presentation.

\begin{theorem}[Semi-sparse $\cK$-Putinar Positivstellensatz]
\label{thm:semiSparse}
Let $\cK$ be a class of polynomials such that $\KSOMS \subseteq \cK$. Let polynomials $g_1,\dots,g_m \in \R[x]$, $r > 0$, and $U_j > 0$ such that $U_j \ge \max\{g_j(x): x \in \B_r\}$, $j=1,\dots,m$, be given. Then, $p(x) > 0$ for all $x \in \{x \in \B_{r}: g_j(x) \ge 0, j=1,\dots,m\}$ implies that there exist $2n+1$-variate polynomial $\alpha_j(y,z,u) \in \cK$ for $j=0,1$, and univariate polynomials $\rho_j(u) \in \cK$, $j=1,\dots,m$, such that
\begin{equation}
\label{eq:semiPutinar}
    p(x) = \alpha_0(x+r e,r e-x,r^2-\|x\|^2) + \alpha_1(x+r e,r e-x,r^2-\|x\|^2)(r^2-\|x\|^2) + \sum_{j=1}^m\rho_j(U_j-g_j(x))g_j(x).
\end{equation}
\end{theorem}

If in Theorem~\ref{thm:semiSparse} one makes the (additional) assumption that the class~$\cK$ contains only non-negative polynomials, then the right hand side of~\eqref{eq:semiPutinar} is clearly a Positivstellensatz for~$p$ on $S=\{x \in \B_{r}: g_j(x) \ge 0, j=1,\dots,m\}$. Weaker assumptions on $\cK$ might also be enough to ensure the non-negativity of the right hand side of~\eqref{eq:semiPutinar} on $S$. For instance, it is enough for $\cK$ to be a class of polynomials that are non-negative on the non-negative orthant (e.g., $\KSAGE$), as $\alpha_0$, $\alpha_1$ and $\rho_j$, $j=1,\dots,m$, are evaluated only on non-negative entries when $x \in S$. We do not need assumptions on~$\cK$ of this type to obtain Theorem~\ref{thm:semiSparse}; that is, even if~\eqref{eq:semiPutinar} is not a Positivstellensatz, for instance when $\cK = \R[x]$, the statement of the theorem is (some times trivially) true.

\begin{myremark}[Sparsity]
\label{rem:sparsity}
Notice that~\eqref{eq:semiPutinar} is written in terms of $\rho_j(U_j - g_j(x))$,  $j=1,\dots,m$, where
$\rho_j$ is a univariate polynomial, making the multivariate polynomial $\rho_j(U_j - g_j(x))$ sparse, as it is constructed from the (few) coefficients of $\rho_j$. On the other hand, the two first terms $\alpha_1$ and $\alpha_2$ in~\eqref{eq:semiPutinar}  are non-sparse. For this reason, we refer to~\eqref{eq:semiPutinar}  as {\em semi-sparse} \germans{}
(see Section~\ref{sec:fullysparse} for more details). Importantly, this term-sparsity is inherent; that is, it arises independent of the sparsity characteristics of the polynomial $p$ and the polynomials $g \in \R[x]^m$ defining the semialgebraic set of interest.
\end{myremark}

We now derive semi-sparse  Putinar-type Positivstellens\"atze for the non-SOS classes of polynomials $\KSONC$, $\KSDSOS$, $\KDSOS$, and $\KSAGE$. These Putinar-type Positivstellens\"atze follow from Theorem~\ref{thm:semiSparse}, as all these classes of polynomials contain the class of SOMS polynomials (Proposition~\ref{prop:relationKs}).
 Moreover, when working with \SAGEtext{} polynomials the Positivstellensatz can be simplified.
\begin{corollary}[Semi-sparse SONC, SDSOS, DSOS and \SAGEtext{} Putinar-type Positivstellensatz]
\label{cor:semiSparse}
Let polynomials $g_1,\dots,g_m \in \R[x]$, $r > 0$, and $U_j > 0$ such that $U_j \ge \max\{g_j(x): x \in \B_r\}$, $j=1,\dots,m$, be given. Let $p(x) > 0$ for all $x \in \{x \in \B_{r}: g_j(x) \ge 0, j=1,\dots,m\}$. Then
\begin{enumerate}[label=(\roman*)]
\item Taking $\cK = \KSONC$,
 $\cK = \KSDSOS$, $\cK = \KDSOS$, or  $\cK = \KSAGE$,
implies that there exist $2n+1$-variate polynomials $\alpha_j(y,z,u) \in \cK$ for $j=0,1$, and univariate polynomials $\rho_j(u) \in \cK$, for $j=1,\dots,m$, such that~\eqref{eq:semiPutinar} holds.
\label{cor:semiSparseNonSAGE}
\item Taking $\cK = \KSAGE$
implies that there exist $\alpha(y,z,u) \in \KSAGE$ and univariate \SAGEtext{} polynomials $\rho_j(u)$ for $j=1,\dots,m$, such that
\begin{equation}\label{eq:PutSAGE}
    p(x) = \alpha(x+re, re-x,r^2-\|x\|^2)+ \sum_{j=1}^m\rho_j(U_j - g_j(x))g_j(x).
\end{equation}
\label{cor:semiSparseSAGE}
\end{enumerate}
\end{corollary}
\begin{proof}
From Proposition~\ref{prop:relationKs}, all the mentioned classes of polynomials contain $\KSOMS$; therefore, Theorem~\ref{thm:semiSparse} implies that~\eqref{eq:semiPutinar} holds, and in particular, statement~\ref{cor:semiSparseNonSAGE} holds. For the case of $\cK = \KSAGE$, statement~\ref{cor:semiSparseSAGE} is obtained from~\ref{cor:semiSparseNonSAGE}, by noticing that $\KSAGE$ is closed under sums and multiplication by a variable. Thus, given that~\eqref{eq:semiPutinar} holds with $\alpha_j(y,z,u) \in \KSAGE$ for $j=0,1$, we have that
$\alpha(y,z,u) := \alpha_0(y,z,u) + \alpha_1(y,z,u)\cdot u \in \KSAGE$, therefore~\eqref{eq:PutSAGE} follows.
\end{proof}

\begin{myremark}[Indeed \germans{}]
\label{rem:RHS2}
it is clear that
Corollary~\ref{cor:semiSparse}\ref{cor:semiSparseNonSAGE} results in \germans{} for $p$ on $S=\{x \in \B_{r}: g_j(x) \ge 0, j=1,\dots,m\}$ for the classes of polynomials $\KSONC$, $\KSDSOS$ and $\KDSOS$, as these are classes of non-negative polynomials. In the case of \SAGEtext{} polynomials, Corollary~\ref{cor:semiSparse}.\ref{cor:semiSparseNonSAGE} (respectively Corollary~\ref{cor:semiSparse}.\ref{cor:semiSparseSAGE}) also results in a \german{} for $p$ on $S$, as the \SAGEtext{} polynomials in~\eqref{eq:semiPutinar} (resp.~\eqref{eq:PutSAGE}) are evaluated on non-negative entries when $x \in S$.
\end{myremark}

Note for example that if the DSOS Putinar-type \german{} of Corollary~\ref{cor:semiSparse} is used in the approximation hierarchies for polynomial optimization problems considered in~\citep[see, e.g.,][]{ahmadi2019dsos, kuang2016alternative, zheng2019sparse}, one would obtain a convergent approximation hierarchy of SOCP problems that (under the assumptions of Corollary~\ref{cor:semiSparse}) are guaranteed to converge to the polynomial optimization problem's objective value.

Recall that one of the motivations for deriving the $\cK$-Putinar \germans{} in Corollary~\ref{cor:semiSparse} for SONC, SDSOS, DSOS and \SAGEtext{} polynomials, is to look for computationally more efficient alternatives to the use of SOS polynomials when constructing  a \german{} for a polynomial over a compact semialgebraic set. In Section~\ref{sec:examples}, we will present examples that illustrate the computational effort required to compute $\cK$-Putinar-type \germans{} based on different classes $\cK$ of polynomials, for some particular instances of polynomials and compact semialgebraic sets.

Next, we look at taking advantage of the sparsity structure of the polynomial $p$ and the polynomials $g\in \R[x]^m$ defining the semialgebraic set $S_g$ to obtain a fully-sparse version of the $\cK$-Putinar \germans{} presented in this section.

\subsection{Fully-sparse non-SOS Putinar-type \germans{}}
\label{sec:fullysparse}
An important contribution in the recent literature related to SOS \germans{}, has been to investigate how to obtain sparse versions of different Positivstellens\"atze. That is, to investigate how to take advantage of the sparsity characteristics of the
the polynomial $p$,  whose non-negativity is being certified, as well as the sparsity of the polynomials $g \in \R[x]^m$ (defining the underlying semialgebraic set $S_g \subseteq \R^n$), to obtain a sparse expression that certifies the non-negativity of the polynomial $p$ on $S$. For example, consider the earlier work in~\citep[to name a few]{parrilo2002explicit, pena2008exploiting, laurent2007semidefinite,waki2006sums}. More recently the focus has been on exploiting {\em correlative and term sparsity} to derive sparse versions of SOS \germans{} such as Putinar's Positivstellensatz~\citep{lasserre2006convergent,wang2020cs,wang2021tssos}; and Reznick's and Putinar-Vasilescu's Positivstellensatz~\citep{mai2020sparse,putinar1999solving}.

Generally, the sparsity of a certificate is classified as either \emph{correlative sparsity} or \emph{term sparsity}. A term is called correlative sparse if it only uses a few variables, and a term is called term sparse if the expression uses only a few terms. For example, box constraints $\{x_i+1\geq 0,1-x_i\geq 0,i=1,\dots,n\}$ are correlative sparse (i.e., each of the terms only depends on one variable), but not term sparse, since $2n$ terms are required to represent the constraints. On the other hand, a ball constraint $r^2-\|x\|^2 \ge 0$, $r > 0$, is not correlative sparse (i.e., each term depends on $n$ variables), but it is clearly term sparse, since only one term is required to represent the constraint.

We can now state the facts in Remark~\ref{rem:sparsity} more formally. The Positivstellensatz of the form~\eqref{eq:semiPutinar}, obtained in Theorem~\ref{thm:semiSparse}, presents term-sparsity
with respect to the terms of the polynomials $g_j$, $j=1,\dots,m$ (as first discussed in Remark~\ref{rem:sparsity}). Namely, $\rho_j$ is a univariate polynomial, and thus $\rho_j(U_j - g_j(x))$ is a term-sparse polynomial for $j=1,\dots,m$. On the other hand, the two first terms $\alpha_1$ and $\alpha_2$ in~\eqref{eq:semiPutinar}  are non-sparse, for this reason we call~\eqref{eq:semiPutinar}  a {\em semi-sparse} certificate. For the case of SOS \germans{}, this type of semi-sparse certificate has been proposed and numerically tested in~\cite{kuryatnikova2019copositive} (see Section~\ref{sec:SOS} for further details).

Under standard correlative structural assumptions on the sparsity of $p$ and the polynomials $g \in \R[x]^m$, we next obtain a {\em fully-sparse} version of the semi-sparse $\cK$-Putinar \germans{} (Theorem~\ref{thm:semiSparse}),
where, besides the term sparsity due to the use of univariate polynomials $\rho_j(U_j-g_j(x))$, the polynomials $\alpha_1$ and $\alpha_2$ in the \germans{} have correlative sparsity as well (see Theorem~\ref{thm:mainSparse}). To formally state this, our second main result, we first need to introduce some relevant definitions related to correlative sparsity~\citep[see, e.g.,][]{lasserre2006convergent}.

For brevity, in what follows, we denote $[n] := \{1,\dots,n\}$. Let $k > 0$,
the family $I=\{I_1,\dots,I_k\}$ is a $k$-cover of $[n]$ if $\bigcup_{l=1}^k I_l = [n]$.
A $k$-cover $I$ of $[n]$ is said to satisfy the running intersection property (RIP) if for all $1 \leq l \leq k-1$ we have $\bigcup_{\ell=1}^{l}I_{\ell}\cap I_{\ell+1}\subseteq I_j$ for some $j\le l$.

Given $x \in \R^n$ and $J \subseteq [n]$, we use $x_J \in \R^J$ to indicate the vector obtained from $x$ by considering only the elements $x_i$ for all $i \in J$. Let $I$ be a given $k$-cover of $[n]$.  A polynomial $p \in \R[x]$ is called $I$-sparse if $p(x)=\sum_{l=1}^kp_l(x_{I_l})$. Also, for $R \in \R^k$,
the set  $\B^I_R = \{x\in \R^n:\|x_{I_l}\| \le R_l,\, l=1,\dots,k\}$ is called the \emph{$I$-sparse ball of radius~$R$}.

With these definitions at hand we can state the fully-sparse $\cK$-Putinar \germans{} , whose proof is delayed to Section~\ref{sec:prooffullysparse} for ease of presentation.

\begin{theorem}[Fully-sparse $\cK$-Putinar Positivstellensatz]\label{thm:mainSparse}
Let $\cK$ be a class of polynomials such that $\KSOMS \subseteq \cK$.
Let $I$ be a $k$-cover of $[n]$ satisfying the RIP. Let $R \in \R_{++}^k$. For $j = 1,\dots, m$, let $l_j \in [k]$,
$g_j \in \R[x_{I_{l_j}}]$,  and $U_j > 0$ such that $U_j \ge \max\{g_j(x_{I_{l_j}}): \|x_{I_{l_j}}\| \le R_{l_j}\}$, be given. Let $p \in \R[x] $ be $I$-sparse. Then $p(x) > 0$ on $\{x \in \B^I_R: g_j(x_{I_{l_j}}) \ge 0, j=1,\dots,m\}$ implies that there exist $(2|I_l|+1)$-variate polynomials $\alpha_{j,l}(y_{I_l},z_{I_l},u) \in \cK $ for $l=1,\dots,k$ and $j=0,1$, and univariate polynomials $\rho_{j}(u)\in\cK$, for  $j= 1,\dots,m$, such that
\begin{multline}
\label{eq:mainSparse}
p(x) = \sum_{l=1}^k \alpha_{0,l}(x_{I_{l}}+R_l e,R_l e - x_{I_{l}},R_l^2 -\|x_{I_{l}}\|^2)
                        + \alpha_{1,l}(x_{I_{l}}+R_l e,R_l e - x_{I_{l}},R_l^2 -\|x_{I_{l}}\|^2)(R_l^2-\|x_{I_{l}}\|^2)\\
                       +\sum_{j=1}^m\rho_{j}(U_j - g_j(x_{I_{l_j}})) g_j(x_{I_{l_j}}).
\end{multline}
\end{theorem}

Similarly to the case of~Theorem~\ref{thm:semiSparse}, the right hand side of~\eqref{eq:mainSparse} in Theorem~\ref{thm:mainSparse} is clearly a \german{} for $p$ on $S = \{x \in \B^I_R: g_j(x_{I_{l_j}}) \ge 0, j=1,\dots,m\}$ when $\cK$ is a class of non-negative polynomials, or more generally a class of polynomials non-negative on the non-negative orthant. Also, similar to Theorem \ref{thm:semiSparse}, in Theorem~\ref{thm:mainSparse}, we have that the terms~$\rho_j(U_j - g_j(x_{I_{l_j}}))$, $j=1,\dots,m$, are term-sparse. In addition, since the polynomials $\alpha_{0,l},\alpha_{1,l}$, $l=1,\dots,k$ in the certificate only depend on part of the variables (namely, on $x_{I_l}$), we have that these polynomials are correlative sparse.

As mentioned earlier, when SOS polynomials are used as the base class of polynomials to certify non-negativity, there are a number of results that show how to exploit correlative sparsity in corresponding SOS Positivstellensatz. In contrast,  the fully-sparse $\cK$-Putinar \germans{} derived above, provide sparse \germans{} in which~$\cK$ is any class of non-negative polynomials, or even certain classes of not-necessarily non-negative polynomials that contain all SOMS polynomials. In particular we obtain sparse Putinar-type \germans{} based on SONC, SDSOS, DSOS and \SAGEtext{} polynomials.
We formally state this fact in the next corollary.

\begin{corollary}[Fully-sparse SONC, SDSOS, DSOS and \SAGEtext{} Putinar-type Positivstellensatz]
\label{cor:mainSparse}
Let $I$ be a $k$-cover of $[n]$ satisfying the RIP. Let $R \in \R_{++}^k$. For $j = 1,\dots, m$, let $l_j \in [k]$,
$g_j \in \R[x_{I_{l_j}}]$,  and $U_j > 0$ such that $U_j \ge \max\{g_j(x_{I_{l_j}}): \|x_{I_{l_j}}\| \le R_{l_j}\}$, be given. Let $p \in \R[x] $ be $I$-sparse and positive on $\{x \in \B^I_R: g_j(x_{I_{l_j}}) \ge 0, j=1,\dots,m\}$.
Then
\begin{enumerate}[label=(\roman*)]
\item Taking $\cK = \KSONC$,
 $\cK = \KSDSOS$, $\cK = \KDSOS$,  or  $\cK = \KSAGE$,
implies that there exist $(2|I_l|+1)$-variate polynomials $\alpha_{j,l}(y_{I_l},z_{I_l},u) \in \cK $ for $l=1,\dots,k$ and $j=0,1$, and univariate polynomials $\rho_{j}(u)\in\cK$, for  $j= 1,\dots,m$, such that~\eqref{eq:mainSparse} holds.
\label{cor:mainSparseNonSAGE}
\item Taking $\cK = \KSAGE$
implies that there exist $\alpha_l \in \KSAGE[y_{I_l},z_{I_l},u]$ for $l=1,\dots,k$, and univariate \SAGEtext{} polynomials $\rho_j$ for $j=1,\dots,m$, such that
\begin{equation}\label{eq:fullyPutSAGE}
    p(x) = \sum_{l=1}^k \alpha_{l}(x_{I_{l}}+R_l e,R_l e - x_{I_{l}},R_l^2 -\|x_{I_{l}}\|^2) +
                       \sum_{j=1}^m\rho_{j}(U_j - g_j(x_{I_{l_j}}))g_j(x_{I_{l_j}}).
\end{equation}
\label{cor:mainSparseSAGE}
\end{enumerate}
\end{corollary}
\begin{proof}
The proof is analogous to that of Corollary~\ref{cor:semiSparse}, but this time using Theorem~\ref{thm:mainSparse}. Namely,
From Proposition~\ref{prop:relationKs}, all the mentioned classes are univariate rich and contain $\KSOMS$; therefore, Theorem~\ref{thm:mainSparse} implies that~\eqref{eq:mainSparse} holds, and in particular, statement~\ref{cor:semiSparseNonSAGE} holds. For the case of $\cK = \KSAGE$, statement~\ref{cor:semiSparseSAGE} is obtained from~\ref{cor:mainSparseNonSAGE}, by noticing that $\KSAGE$ is closed under sums and multiplication by a variable. Thus, given that~\eqref{eq:mainSparse}  holds with $\alpha_{j,l}(y_{I_l},z_{I_l},u) \in \KSAGE$ for $l=1,\dots,k$ and $j=0,1$, and univariate polynomials $\rho_j(u) \in \KSAGE$, for $j=1,\dots,m$, we have that
$\alpha_l(y_{I_l},z_{I_l},u) := \alpha_{0,l}(y_{I_l},z_{I_l},u) + \alpha_{1,l}(y_{I_l},z_{I_l},u)\cdot u \in \KSAGE$ and $\tau_j(u) := \rho_j(u)\cdot u \in \KSAGE$, and~\eqref{eq:fullyPutSAGE} follows.
\end{proof}

Clearly, Remark~\ref{rem:RHS2} applies to Corollary~\ref{cor:mainSparse}, that is we do obtain a fully-sparse \german{} for each of the SONC, SDSOS, DSOS, SOMS and \SAGEtext{} polynomial classes.

\subsection{Illustrative Examples}
\label{sec:examples}

We used Example~\ref{ex:counterSDSOS} to illustrate the fact that a non-SOS Putinar-type \german{} cannot be obtained simply by replacing the SOS polynomials in~\eqref{eq:putinarform} by SDSOS polynomials. This is a strategy that was proposed for numerical purposes in~\citep{ahmadi2019dsos} to approximately solve polynomial optimization problems with lower computational effort than by using the well-known approach of~\citet{lasserre2001global}.  Corollary~\ref{cor:semiSparse}\ref{cor:semiSparseNonSAGE} provides a (semi-sparse) SDSOS Putinar-type \german{} of the general form~\eqref{eq:putinarform}, in which the class of SDSOS polynomials is used as the base class to certify non-negativity. To illustrate this, in Example~\ref{ex:counterSDSOSrev}, we show that the non-negativity of the polynomial over the semialgebraic set considered in Example~\ref{ex:counterSDSOS} can be evidenced using SDSOS polynomials.

\begin{myexample}[Example~\ref{ex:counterSDSOS} revisited: SOMS Putinary-type \german{}]
\label{ex:counterSDSOSrev}
Let $g(x_1,x_2) = 1-x_1^2-x_2^2$.
Recall from Example~\ref{ex:counterSDSOS} that the ``simple" DSOS/SDSOS Putinar Positivstellensatz fails to certify non-negativity of the polynomial
$p(x_1,x_2)=(2-x_1-x_2)^2$ on $S_g = \{ (x_1,x_2) \in \R^2: 1-x_1^2-x_2^2 \ge 0\} = \B_{1}$ even though $p > 0$ on $\B_1$. Let $m= 0$. By Theorem~\ref{thm:semiSparse},
there exist $\alpha_0(y_1,y_2,z_1,z_2,u), \alpha_1(y_1,y_2,z_1,z_2,u) \in \KSOMS \subset \KSDSOS$ such that
\begin{multline}
\label{eq:ex3}
(2-x_1-x_2)^2=\alpha_0(x_1+1,x_2+1,1-x_1,1-x_2,1-x_1^2-x_2^2)\\
+ \alpha_1(x_1+1,x_2+1,1-x_1,1-x_2,1-x_1^2-x_2^2)(1-x_1^2-x_2^2).
\end{multline}
Indeed, it is easy to check that one can choose
\[
\begin{array}{ll}
\alpha_0(y_1,y_2,z_1,z_2,u)=&\tfrac{1}{4}+\tfrac{1}{2}z_2^2+\tfrac{1}{4}u^2+\tfrac{1}{3}z_1^4
+\tfrac{7}{12}(z_1z_2)^2+\tfrac{1}{12}(z_1y_1)^2
+\tfrac{1}{12}(z_1y_2)^2+\tfrac{1}{4}z_2^4\\[2ex]
\alpha_1(y_1,y_2,z_1,z_2,u)=&\tfrac{1}{2}+\tfrac{2}{3}z_1^2+\tfrac{1}{2}z_2^2.
\end{array}
\]
Note that the right hand side of~\eqref{eq:ex3} makes the non-negativity of $p$ on $S_g$ evident, since $\alpha_0$ and $\alpha_1$ are non-negative everywhere, and $1-x_1^2-x_2^2 \ge 0$ for all $x \in S_g$.
\end{myexample}

As mentioned earlier, Corollary~\ref{cor:mainSparse} shows that the presence of term and correlative sparsity can be used to obtain sparse \germans{} that are not necessarily obtained using SOS polynomials as the base class to certify non-negativity. Example~\ref{ex:sparseSONC} illustrates this fact by showing how a sparse \german{}, using SDSOS polynomials, can be constructed by taking advantage of term and correlative sparsity.

\begin{myexample}[Fully-sparse SDSOS Putinar-type \german{}]
\label{ex:sparseSONC}
Let $n=3$, $g_i(x) = 1-x_i$, $g_{n+i}(x) = 1+x_i$ for $i=1,\dots,n$, and $p(x_1,x_2,x_3)=11+2x_1^2+4x_1x_2-x_2^2-2x_2x_3-3x_3-2x_3^3$.
Note that for any $(x_1, x_2, x_3) \in S_g=[-1,1]^3$, $p(x_1,x_2,x_3) \ge 11+2x_1^2-4x_1-1-2-3-2 \ge 1 >0$.
Also, both the polynomial $p$ and the polynomials $g \in \R[x]^{2n}$ generating the set $S_g$ are $I$-sparse with $I=\{\{1,2\},\{2,3\}\}$, which satisfies the RIP (c.f., Section~\ref{sec:fullysparse}). Further, notice that by setting $R =(\sqrt{2},\sqrt{2})\tr$, it follows that $S_g \cap \B^I_R = S_g$, where
$\B^I_R = \{(x_1,x_2,x_3) \in \R^3: \|(x_1,x_2)\| \le \sqrt{2}, \|(x_2,x_3)\|\le \sqrt{2}\}$. Let $U_j=3\geq\max\{g_j(x):x\in B^I_R\}$ for $j=1,\dots,2n$. Then by
Corollary~\ref{cor:mainSparse}\ref{cor:mainSparseNonSAGE}, there exists
$\alpha_{j,l} \in \KSDSOS $ for $l=1,2$ and $j=0,1$, and univariate polynomials $\rho_{i,1}, \rho_{i,2} \in\KSDSOS$, for  $i= 1,\dots,3$, such that
\begin{equation}
\label{eq:SONCcert}
\begin{array}{lcl}
    p(x)&=&\alpha_{0,1}(x_1 + \sqrt{2}, x_2 + \sqrt{2}, \sqrt{2}-x_1,\sqrt{2}-x_2,2-x_1^2-x_2^2)  \\
    &&\quad+ \alpha_{1,1}(x_1 + \sqrt{2}, x_2 + \sqrt{2}, \sqrt{2}-x_1,\sqrt{2}-x_2,2-x_1^2-x_2^2)(2-x_1^2-x_2^2) \\
    &&\quad\quad+ \alpha_{0,2}(x_2 + \sqrt{2}, x_3 + \sqrt{2}, \sqrt{2}-x_2,\sqrt{2}-x_3,2-x_2^2-x_3^2)\\
    &&\quad\quad\quad+ \alpha_{1,2}(x_2 + \sqrt{2}, x_3 + \sqrt{2}, \sqrt{2}-x_2,\sqrt{2}-x_3,2-x_2^2-x_3^2)(2-x_2^2-x_3^2)\\
    &&\quad\quad\quad\quad+\sum_{i=1}^3 \left (\rho_{i,1}(2+x_i)(1-x_i) + \rho_{i,2}(2-x_i)(1+x_i) \right).
    \end{array}
    \end{equation}
Indeed, it is easy to check that one can choose
\[
\begin{array}{rcl@{\hskip 0.5in}rcl}
    \alpha_{0,1}(y_1,y_2,z_1,z_2,u)&=&\frac{1}{4}(y_1y_2-z_1z_2)^2&    \rho_{1,1}(u)&=&0\\
    \alpha_{1,1}(y_1,y_2,z_1,z_2,u)&=&0&    \rho_{1,2}(u)&=&0\\
    \alpha_{0,2}(y_2,y_3,z_2,z_3,u)&=&\frac{1}{8}(y_2z_3-z_2y_3)^2&    \rho_{2,1}(u)&=&0.5(2-u)^2\\
    \alpha_{1,2}(y_2,y_3,z_2,z_3,u)&=&5&    \rho_{2,2}(u)&=&0.5(2-u)^2\\
    &&&\rho_{3,1}(u)&=&2(2.5-u)^2+0.5\\
    &&&\rho_{3,2}(u)&=&0.
\end{array}
\]
Note that since
$\alpha_{j,l}$, $l=1,2$ and $j=0,1$ and $\rho_{i,1}, \rho_{i,2}$, for  $i= 1,\dots,3$ are non-negative polynomials, and the polynomials $2-x_1^2-x_2^2, 2-x_2^2-x_3^2, 1- x_i,1+x_i$ are non-negative in $S_g$, expression~\eqref{eq:SONCcert} evidences the non-negativity of $p$ on $S_g$.
\end{myexample}

There is a clear computational motivation behind the development of
semi-sparse (see, e.g.,~\citep[Cor. 7]{kuryatnikova2019copositive} and Proposition~\ref{prop:KputSOS})
and fully-sparse~\citep[see, e.g.,][Cor. 3.3]{lasserre2006convergent} SOS \germans{} that take advantage of term or both term and correlative sparsity. That is, loosely speaking, the SDP that needs to be solved to construct a sparse SOS \german{} is substantially smaller than the one that needs to be solved to construct the corresponding non-sparse SOS \german{}. Corollaries~\ref{cor:semiSparse} and~\ref{cor:mainSparse} show that semi-sparse and fully-sparse Putinar-type \germans{} can be constructed using non-SOS classes of polynomials such as SONC, SDSOS, DSOS, and \SAGEtext{} polynomials. Such certificates substantially reduce the complexity of the optimization problem (e.g., a SOCP problem for a SDSOS \german{}, or a REP problem for a \SAGEtext{} \german{}) that needs to be solved to construct the corresponding \german{}. In Example~\ref{ex:sparseSAG}, we illustrate this for the case in which the goal is to construct SONC, SDSOS, DSOS and/or \SAGEtext{} Putinary-type \germans{} for the generalized Rosenbrock function over the hypercube.

\begin{myexample}[Sparsity vs size of the \germans{}]
\label{ex:sparseSAG}
Let $g_i(x) = 1-x_i$, and $g_{n+i}(x) = 1+x_i$ for $i=1,\dots,n$. Consider the generalized Rosenbrock function
\[
    p(x)=\sum_{i=1}^{n-1}[100(x_{i+1}-x_i^2)^2+(1-x_{i+1})^2],
\]
and the problem of certifying its non-negativity on the hypercube $S_g=[-1,1]^n$. The polynomial $p(x)$  and the polynomials defining the hypercube are $I$-sparse with respect to $I=\{I_1,\dots,I_{n-1}\}$, where $I_i=\{i,i+1\}$ for $i=1,\dots,n-1$.
Note that $I$ satisfies the RIP (c.f., Section~\ref{sec:fullysparse}). Further, notice that by setting $R =(\sqrt{2},\dots,\sqrt{2})\tr \in \R^{n-1}$, it follows that $S_g  \subseteq \B^I_R = \{x\in \R^n: \|(x_i,x_{i+1})\| \le \sqrt{2}, i=1,\dots,n-1\}$. Also, similar to Example~\ref{ex:sparseSONC}, let $U_j = 3$, for $j=1,\dots,2n$. Then by
Corollary~\ref{cor:mainSparse}\ref{cor:semiSparseNonSAGE}, there exist fully-sparse $\cK$-Putinar \germans{} for $p$ over $S_g$ for $\cK \in \{\KSONC, \KSDSOS, \KDSOS, \KSAGE\}$ of the form~\eqref{eq:mainSparse}. Also, if we ignore the sparsity structure defined by $I$,
by Corollary~\ref{cor:semiSparse}\ref{cor:mainSparseNonSAGE}, there exist semi-sparse $\cK$-Putinar \germans{} for $p$ over $S_g$ for $\cK \in \{\KSONC, \KSDSOS,  \KDSOS, \KSAGE\}$ of the form~\eqref{eq:semiPutinar}. In Table~\ref{table:SONC-GRos}, we compare the number and size of the polynomials, in terms of both the dimension of their corresponding Gram matrices (1st row), and the number of monomials (2nd row), that are required to test for the existence of the \germans{}, when variable polynomials in the
 \german{} expressions (e.g., $\alpha_0$, $\alpha_1$,  and $\rho_j$, $j=1,\dots,m$ in~\eqref{eq:semiPutinar}) of degree up to $2d$ are considered. In the fully-sparse (resp. semi-sparse) case; that is, the 3rd column (resp. 2nd column) in Table~\ref{table:SONC-GRos}, the entries on the table correspond to the number and size of the polynomials needed to construct expression~\eqref{eq:mainSparse} (resp. expression~\eqref{eq:semiPutinar}). In the non-sparse case; that is, the 1st column in Table~\ref{table:SONC-GRos}, the entries on the table correspond to the number and size of the polynomials needed to construct expression~\eqref{eq:semiPutinar}, {\em \bf but} without the assumption of the polynomials $\rho_j$, $j=1,\dots,m$ in~\eqref{eq:semiPutinar} being univariate.
For the particular case of $n=4$, $d=2$, the corresponding size values are given numerically in the Table. Furthermore, in this particular case the fully-sparse SONC, SDSOS, and \SAGEtext{} \germans{} expressions for $p(x) + \epsilon$;
that is, the coefficients of the corresponding variable polynomials, were computed with $\epsilon = 0.001$ using SOCP, SOCP, and REP respectively. These coefficients are not presented for the purpose of brevity.

\begin{table}[!htb]
\centering
\resizebox{\textwidth}{!}{%
\begin{tabular}{lccccc}
\toprule
& && \multicolumn{3}{c}{Non-SOS \germans{}} \\
\cmidrule{4-6}
                  &$n$ & $d$  & Non-sparse          & Semi-sparse~\eqref{eq:semiPutinar} & Fully-sparse~\eqref{eq:mainSparse}\\
 \midrule
 \multirow{2}{*}{Size of Gram matrices}    & $\cdot$ & $\cdot$ & $2{2n+1+d \choose 2n+1}^2+2n{n+1+d \choose n+1}^2$  & $2{2n+1+d \choose 2n+1}^2+2n(1+d)^2$       &    $2(n-1){5+d \choose 5}^2 +2n(1+d)^2$     \\
& $4$ & $2$   &= 7850         &= 6122       &   = 2718     \\
 \midrule
 \multirow{2}{*}{Number of monomials}   & $\cdot$ & $\cdot$ & $2{2n+1+2d \choose 2n+1}+2n{n+1+2d \choose n+1}$  & $2{2n+1+2d \choose 2n+1}+2n(1+2d)$       &    $2(n-1){5+2d \choose 5} + 2n(1+2d)$     \\
 & $4$ & $2$ &= 1990         &= 1470       &= 796     \\
 \bottomrule
\end{tabular}
}
\caption{Size of non-SOS Putinar-type \germans{} expressions with variable polynomials up to degree $2d$ for the generalized Rosenbrock function.}
\label{table:SONC-GRos}
\end{table}

Table~\ref{table:SAGE-GRos} is similar to Table~\ref{table:SONC-GRos}, but now for the case of \SAGEtext{} Putinar-type \germans{}, in which advantage is taken of the fact that \SAGEtext{} polynomials are closed under multiplication by monomials (see Corollary~\ref{cor:semiSparse}\ref{cor:semiSparseSAGE} and Corollary~\ref{cor:mainSparse}\ref{cor:mainSparseSAGE}). That is, the entries in the fully-sparse (resp. semi-sparse) column of the table now refer to expression~\eqref{eq:fullyPutSAGE} (resp. expression~\eqref{eq:PutSAGE}), and the entries in the non-sparse column refer to expression~\eqref{eq:PutSAGE}, {\em \bf but} without the assumption of the polynomials $\rho_j$, $j=1,\dots,m$ in~\eqref{eq:semiPutinar} being univariate. For the particular case of $n=4$, $d=2$, the corresponding size values are given numerically in the Table. Furthermore, in this particular case the the fully-sparse \SAGEtext{} \german{} expression  for $p(x)+ \epsilon$;
that is, the coefficients of the corresponding variable polynomials, were computed up to numerical precision using REP with $\epsilon = 0.001$. These coefficients are not presented for the purpose of brevity.

\begin{table}[!htb]
\centering
\resizebox{\textwidth}{!}{%
\begin{tabular}{lccccc}
\toprule
& && \multicolumn{3}{c}{\SAGEtext{} \germans{}} \\
\cmidrule{4-6}
                  &$n$ & $d$  & Non-sparse          & Semi-sparse~\eqref{eq:PutSAGE} & Fully-sparse~\eqref{eq:fullyPutSAGE}\\
 \midrule
 \multirow{2}{*}{Size of Gram matrices}   & $\cdot$ & $\cdot$ & ${2n+1+d \choose 2n+1}^2+2n{n+1+d \choose n+1}^2$  & $2{2n+1+d \choose 2n+1}^2+2n(1+d)^2$       &    $(n-1){5+d \choose 5}^2+2n(1+d)^2$     \\
 & $4$ & $2$ &= 4825         &= 3097       &   = 1395     \\
 \hline
\multirow{2}{*}{Number of monomials}    & $\cdot$ & $\cdot$ & ${2n+1+2d \choose 2n+1}+2n{n+1+2d \choose n+1}$  & $2{2n+1+2d \choose 2n+1}+2n(1+2d)$       &    $(n-1){5+2d \choose5}+2n(1+2d)$     \\
& $4$ & $2$ &= 1275         &= 755       &   = 418     \\
\bottomrule
\end{tabular}
}
\caption{Size of \SAGEtext{} Putinary-type \germans{} expressions with variable polynomials up to degree $2d$ for the generalized Rosenbrock function.}
\label{table:SAGE-GRos}
\end{table}
\end{myexample}

Three main conclusions can be derived from Example~\ref{ex:sparseSAG}. First, that even if the polynomial $p$ and the polynomials $g \in \R[x]^m$ defining the semilagebraic set $S_g$ do not exhibit sparsity, it is advantageous to use the semi-sparse $\cK$-Putinar \germans{} provided by Corollary~\ref{cor:semiSparse}, in which some of the terms of the \germans{} are constructed using univariate polynomials. Second, that
the fully-sparse $\cK$-Putinar \germans{} provided by Corollary~\ref{cor:mainSparse} allows to take advantage of term and correlative sparsity exhibited by the polynomial $p$ and the polynomials $g \in \R[x]^m$ defining the semilagebraic set $S_g$, to obtain significant reductions in the size of the \germans{},  in the case when the class of polynomials using as the base to certify non-negativity are the (non-SOS) classes of SONC, SDSOS, DSOS, and \SAGEtext{} polynomials. Third, in terms of the size of the polynomials involved in the \germans{}, the fact that \SAGEtext{} polynomials are closed under multiplication by monomials, is advantageous (compare Table~\ref{table:SONC-GRos} and Table~\ref{table:SAGE-GRos}), as the number of terms in the \german{} significantly decreases (see Corollaries~\ref{cor:semiSparse}\ref{cor:semiSparseSAGE} and~\ref{cor:mainSparse}\ref{cor:mainSparseSAGE}).

\section{Delayed proof details and the case of SOS polynomials}
\label{sec:proofs}

The proof of our main results follow the ideas in~\citep{kuryatnikova2019copositive,pena2008exploiting} of reducing questions about the existence of Positivstellens\"atze on semialgebraic sets, to questions about the existence of \germans{} on ``simpler'' sets. Next, in Sections~\ref{sec:proofsemisparse} and~\ref{sec:prooffullysparse}, we use these ideas to present the details of the proofs that were delayed in previous sections. Further, given that SOS polynomials clearly contain all SOMS polynomials, in Section~\ref{sec:SOS}, we consider the results obtained after applying both the semi-sparse and fully-sparse $\cK$-Putinar \germans{} presented in Section~\ref{sec:about} (i.e., Theorems~\ref{thm:semiSparse} and~\ref{thm:mainSparse}) to the special case in which $\cK = \KSOS$.

\subsection{Proof of semi-sparse $\cK$-Putinar \germans{} (Theorem~\ref{thm:semiSparse})}
\label{sec:proofsemisparse}

To prove Theorem~\ref{thm:semiSparse}, in Lemma~\ref{lem:StoB}, we first show that given a compact set $C$ and polynomials $g_j \in \R[x]$, $j=1,\dots,m$, the problem of constructing \germans{} for a polynomial over the set $\{x \in C: g(x) \ge 0\}$ could be reduced to constructing a \german{} for an appropriate polynomial over~$C$.
We then obtain the semi-sparse $\cK$-Putinar \germans{} (Theorem~\ref{thm:semiSparse}) by applying Lemma~\ref{lem:StoB} with $C = \B_r$, after proving the result for the special case in which the semialgebraic set is the ball $\B_r$. That is, Lemma~\ref{lem:StoB} reduces all the effort of the proof to the ball case. Indeed, in Lemma~\ref{lem:KputBall}, the result for the special ball case is proven  using the semi-sparse $\cK$-\schmu{} \german{} (Theorem~\ref{Thm:Kschmudgen}). This proof outline, as well as the formal statements and proofs of Lemma~\ref{lem:StoB} and Lemma~\ref{lem:KputBall} are stated next.

\begin{lemma}
\label{lem:StoB}
Let $\cK$ be a class of polynomials such that $\KSOMS \subseteq \cK$. Let a compact set
$C \subset \R^n$, polynomials $g_1,\dots,g_m \in \R[x]$, and
$U_j > 0$ such that $U_j \ge \max\{g_j(x): x \in C\}$, $j=1,\dots,m$, be given.  Then $p(x) > 0$ on $\{x \in C: g_j(x) \ge 0, j=1,\dots,m\}$ implies that there exist univariate polynomials $\rho_j(u) \in \cK$, $j=1,\dots,m$, such that
\[
p(x) - \sum_{j=1}^m\rho_j(U_j-g_j(x))g_j(x) > 0 \text{ for all $x \in C$}.
\]
\end{lemma}
\begin{proof}Let $S = \{x \in C: g_j(x) \ge 0,\, j=1,\dots,m\}$ and let $\epsilon>0$ be such that $p(x)>\epsilon$ for all $x\in S$.
 Let $M = \min\{p(x):x \in C\}$.
If $M > 0$, the statement follows by letting $\rho_j(u) = 0$ for $j = 1,\dots,m$. Thus, we assume that $M \le 0$.

Now,  let $\delta := \max\{\min\{g_j(x): j=1,\dots,m\}: x\in C,\, p(x) \le \epsilon/2\}$.
We claim $\delta<0$. To prove this claim by contradiction, assume $\delta\geq 0$. The function $\min\{g_j(x):j=1,\dots,m\}$ is continuous and the set $X=\{x\in C:p(x)\leq \epsilon/2\}$ is compact and non-empty. Therefore the value $\delta$ is attained; that is, there is $x^*\in X$ such that $\min\{g_j(x^*):j=1,\dots,m\}=\delta\geq 0$. Thus $x^*\in S$, and then $p(x^*)>\epsilon$ contradicting $x\in X$.

Now, for each $j=1,\dots,m$, let
\[
d_j \in \mathbb{N} \text{ such that }   d_j \ge \hat{d}_j := \frac{\ln(4(m+1)U_j(-M+\epsilon))-\ln(|\delta|\epsilon)}{2\ln(1+|\delta|/U_j)} \text{ and } c_j=\frac{\epsilon}{4(m+1)U_j^{2d_j+1}} > 0.
\]
Define $\rho_j(u)=c_ju^{2d_j}\in \KSOMS\subseteq\cK$. For $t\in[0,U_j]$ we have $\rho_j(U_j-t)\leq c_jU_j^{2d_j}= \frac{\epsilon}{4(m+1)U_j}$. Also, for $t\leq \delta$ we have
 \[
 \rho_j(U_j-t) \ge  c_j(U_j+|\delta|)^{2d_j} = \frac{\epsilon(U_j+|\delta|)^{2d_j}}{4(m+1)U_j^{2d_j+1}}  \ge
 \frac{ \epsilon(1+|\delta|/U_j)^{2\hat{d}_j}}{4(m+1)U_j} = \frac{-M+\epsilon}{|\delta|}.
 \]
We claim that for any $x \in C$, $p(x)-\sum_{j=1}^m \rho_j(U_j-g_j(x))g_j(x) \geq \frac {\epsilon}{4}$. To show this take $x\in C$. Notice  that if $p(x) \ge \epsilon/2$ we have  $p(x)-\sum_{j=1}^m \rho_j(U_j-g_j(x))g_j(x) \ge \frac{\epsilon}{2} - \sum_{j=1}^m \frac{\epsilon}{4(m+1)} > \frac {\epsilon}4$.
If $p(x) \le \epsilon/2$, there is $j^* \in \{1,\dots,m\}$ such that $g_{j^*}(x) \le \delta$ and then $-\rho_{j^*}(U_{j^*}-g_{j^*}(x))g_{j^*}(x) \ge \frac{-M + \epsilon}{-|\delta|}\delta  = -M + \epsilon$.
Therefore $p(x)-\sum_{j=1}^m \rho_j(U_j-g_j(x))g_j(x) \ge M +(-M + \epsilon) - \sum_{j\neq j^*}^m \frac{\epsilon}{4(m+1)U_j}U_j > \epsilon/4$.
\end{proof}

Notice that if one obtains a non-SOS Putinar-type Positivstellensatz for positive polynomials over the ball, then Lemma~\ref{lem:StoB} can be used to extend such result to any compact semialgebraic set. Next we obtain the former result, using the $\cK$-Schm\"{u}dgen \germans{}~\citep[][Cor.~5]{kuryatnikova2019copositive} (Theorem~\ref{Thm:Kschmudgen}).

\begin{lemma}[$\cK$-Putinar \germans{} over the ball]\label{lem:KputBall}
Let $\cK$ be a class of non-negative polynomials satisfying that $\KSOMS \subseteq \cK$.
Let $r >0 $. Then $p(x) > 0$ on $\B_{r}$ implies that there exist $2n+1$-variate polynomials $\alpha_j(y,z,u) \in \cK$ for $j=0,1$ such that
\begin{equation*}
    p(x) = \alpha_0(x+r e,r e-x,r^2-\|x\|^2) + \alpha_1(x+r e,r e-x,r^2-\|x\|^2)(r^2-\|x\|^2).
\end{equation*}
\end{lemma}
\begin{proof}
From Theorem~\ref{Thm:Kschmudgen} applied to $\cK = \R^n_+$, if $p(x)>0$ for all $x\in \B_r$, then there exist $d \ge 0$ and $c_{\alpha,\beta,\gamma} \ge 0$ for all $(\alpha,\beta,\gamma) \in \N^{2n+1}_d$ such that
\begin{equation*}
    p(x) = \sum_{(\alpha,\beta,\gamma) \in \N^{2n+1}_d} c_{\alpha,\beta,\gamma} \left (x+\frac{n+1}2 r e \right )^\alpha \left (\frac{n+1}2 r e-x \right )^\beta (r^2-\|x\|^2)^\gamma.
\end{equation*}
Using $2r (\frac{n+1}2r \pm x_i) = (r \pm x_i)^2 + (r^2 - \|x\|^2)+\sum_{j \neq i} \frac 12 ((r - x_j)^2 + (r + x_j)^2)$
we obtain an expression of degree $2d$ of the form
\begin{equation*}
\begin{split}
    p(x) = \sum_{(\alpha,\beta,\gamma) \in \N^{2n+1}_d} c'_{\alpha,\beta,\gamma}(x+r e)^{2\alpha} \left (r e-x \right)^{2\beta} (r^2-\|x\|^2)^{2\gamma} + \\\sum_{(\alpha,\beta,\gamma) \in \N^{2n+1}_d} c''_{\alpha,\beta,\gamma}(x+r e)^{2\alpha} \left (r e-x \right)^{2\beta} (r^2-\|x\|^2)^{2\gamma+1},
    \end{split}
\end{equation*}
where all the coefficients are non-negative.
From the fact that $\KSOMS \subseteq \cK$, it follows that $y^{2\alpha}z^{2\beta} u^{2\gamma} \in \cK$ for each $(\alpha,\beta,\gamma)$ and the statement follows by collecting terms.
\end{proof}

Using Lemma~\ref{lem:StoB} to extend the $\cK$-Putinar \germans{} over the ball (Lemma~\ref{lem:KputBall}), the semi-sparse $\cK$-Putinar \germans{} (Theorem~\ref{thm:semiSparse}) follows.
\begin{proof}[of semi-sparse $\cK$-Putinar \germans{} (Theorem~\ref{thm:semiSparse})]
Applying Lemma \ref{lem:StoB}, for the case $C = \B_r$; that is,
to reduce the result to certifying positivity over the ball $ \B_r$,
there exist univariate polynomials $\rho_j(u)\in \cK$, $j=1,\dots,m$, such that
\[
p(x)-\sum_{j=1}^m\rho_j(U_j - g_j(x))g_j(x) > 0 \text{ for all $x \in \B_{r}$}.
\]
Then, from Lemma~\ref{lem:KputBall}, there exist $2n+1$-variate polynomials $\alpha_j(y,z,u) \in \cK$ for $j=1,\dots,2$ such that
\[
p(x)-\sum_{j=1}^m\rho_j(U_j - g_j(x))g_j(x) = \alpha_0(x+r e,r e-x,r^2-\|x\|^2) + \alpha_1(x+r e,r e-x,r^2-\|x\|^2)(r^2-\|x\|^2).
\]
\end{proof}

\subsection{Proof of  fully-sparse $\cK$-Putinar \germans{} (Theorem~\ref{thm:mainSparse})}
\label{sec:prooffullysparse}

The semi-sparse $\cK$-Putinar \germans{} (Theorem~\ref{thm:semiSparse}) have inherent term sparsity, regardless of the polynomial $p$ and the polynomials $g \in \R[x]^m$ of interest (see Remark~\ref{rem:sparsity}). Further, when  sparsity information about the polynomial $p$ and the polynomials $g \in \R[x]^m$ of interest is known, we obtain a  fully-sparse $\cK$-Putinar \germans{} (Theorem~\ref{thm:mainSparse}) using a similar approach as the one used to prove the semi-sparse $\cK$-Putinar \germans{}, but considering $C=
\B^I_R$ the $I$-sparse ball. Proposition~\ref{prop:RIPsparse} shows that if $I$ satisfies the RIP property,  we obtain the fully-sparse $\cK$-Putinar \germans{} for the special $I$-sparse ball case.

\begin{proposition}\label{prop:RIPsparse}
Let $I=\{I_1,\dots,I_k\}$ be a $k$-cover of $[n]$ satisfying the RIP. Let $R \in \R^k_{++}$. Let $p(x)$ be an $I$-sparse polinomial positive on $\B^I_R$. Then, for $l=1,\dots,k$, there exist $\hat{p}_l \in \R[x_{I_l}]$ positive on $\B^{I_l}_{R_l}$ such that $p(x)=\sum_{l=1}^k\hat{p}_l(x_{I_l})$.
\end{proposition}
\begin{proof}
The proof is by induction on $k$. For $k=1$ the proposition automatically follows by taking $\hat{p_1}(x)=p(x)$.
Suppose the proposition holds for $k-1$.
Let $\epsilon > 0$ be  such that $p(x)>\epsilon$ for all $x\in \B^I_{R}$. Write $p(x) = \sum_{l=1}^k p_l(x_{I_l})$.
Let $J = I_k\cap \bigcup_{l=1}^{k-1} I_l$.
We write $p_k(x_{I_k})$ as $p_k(x_J,x_{I_k \setminus J})$.
Let
\[
f(x_J)= \min_{y \in \R^{I_k\setminus J}: \|y\|^2 + \|x_J\|^2 \le R_k^2} p_k(x_J,y),
\]
which is a continuous function~\citep[see, e.g.,][]{fiacco1980continuity}.
Using the Stone-Weierstrass Theorem \cite{de1959stone}, there is a polynomial~$q \in \R[x_J]$ such that $|f(x_J) - q(x_J)| \le \epsilon/4 $ for all $x_J$ such that $\|x_J\| \le R_k$. Define $ \hat p_k(x_{I_k})=p_k(x_{I_k})- q(x_J) + \epsilon/2$. For any $x_{I_k} \in \B^{I_k}_{R_k}$ we have
\begin{equation}
\label{eq:hatpk}
\hat p_k(x_{I_k}) = p_k(x_{I_k})- q(x_J) + \epsilon/2 \ge p_k(x_J,x_{I_k-J})- f(x_J) + \epsilon/4 \ge \epsilon/4.
\end{equation}
Let $n_1 = |I_k \setminus J|$. By reordering the indexes of the variables we can assume, without loss of generality, that $I_k \setminus J = \{n-n_1+1,\dots,n\}$; that is, the set of the last $n_1$ indexes. We have then that $\tilde I = \{I_1,\dots,I_{k-1}\}$ is a $(k-1)$-cover of $[n-n_1]$ that satisfies the RIP. Also, from the fact that $I$ satisfies the RIP, we have that $J \subset I_l$ for some $l<k$. Let $\tilde p(x_{[n-n_1]}) = \sum_{l=1}^{k-1}p_l(x_{I_l}) + q(x_J) - \epsilon/2$. Note that $\tilde p$ is $\tilde I$-sparse and  $p(x) = \tilde p(x_{[n-n_1]}) +\hat p_k(x_{I_k})$.

Let $\tilde R = (R_1,\dots,R_{k-1}) \in \R^{k-1}_+$ be obtained by dropping the last entry in $R$. We claim that $\tilde p$ is positive on $\B^{\tilde I}_{\tilde R} \subset \R^{[n-n_1]}$ and thus, by induction hypothesis, $\tilde p(x) = \sum_{l=1}^{k-1}\hat p_l(x_{I_l})$ with $\hat p_l(x_{I_l}) >0$ on $B^{I_l}_{R_l}$ for $l=1,\dots,k-1$. Writing $p(x) = \tilde p(x_{[n-n_1]}) + \hat p_k(x_{I_k})$ and using~\eqref{eq:hatpk}, the statement of the proposition follows.

Now we prove the claim. Because $\B^{\tilde I}_{\tilde R}$ is compact, there is $x^* \in \B^{\tilde I}_{\tilde R}$ such that $\tilde p(x^*) = \min \{ \tilde p(x): x \in  \B^{\tilde I}_{\tilde R}\}$. Also by compactness there is $y^* \in \R^{I_k\setminus J}$ such that $(x^*_J,y^*) \in \B^{I_k}_{R_k}$ and $f(x^*_J) = p_k(x^*_J,y^*)$. We have then that $(x^*,y^*) \in \B^I_R$ and $p(x^*,y^*) = \tilde p(x^*) +\hat p_k(x^*_J,y^*)$. Thus
\[\hat p_k(x^*_J,y^*) = p_k(x^*_J,y^*) - q(x^*_J) + \epsilon/2 = f(x^*_J) - q(x^*_J) + \frac{\epsilon}{2} \le \frac{3}{4}\epsilon,\]
and therefore for any $x \in \B^{\tilde I}_{\tilde R}$,
\[\tilde p(x) \ge \tilde p(x^*) = p(x^*,y^*) - \hat p_k(x^*_J,y^*) \ge \frac{\epsilon}{4}.\]
\end{proof}

With Proposition~\ref{prop:RIPsparse}, the proof of Theorem~\ref{thm:mainSparse} is analogous to that of Theorem~\ref{thm:semiSparse}.

\begin{proof}[of fully-sparse $\cK$-Putinar \germans{} (Theorem~\ref{thm:mainSparse})]
Applying Lemma~\ref{lem:StoB}, for the case $C = \B^I_R$; that is,
to reduce the result to certifying positivity over the ball $\B^I_R$, there exist univariate polynomials $\sigma_j(u)\in \cK$, $j=1,\dots,m$, such that
\begin{equation}
\label{eq:sparse1}
q(x):= p(x)-\sum_{j=1}^m\sigma_j(U_j - g_j(x_{I_{l_j}}))g_j(x_{I_{l_j}}) > 0 \text{ for all $x \in \B^I_R$}.
\end{equation}
Notice that $q(x)$ is $I$-sparse. Thus, from Proposition~\ref{prop:RIPsparse}, it follows that $q(x) = \sum_{l=1}^{k} q_l(x_{I_l})$ where $q_l(x_{I_l}) > 0$ on $\B^{I_l}_{R_l}$.
Then, from Lemma~\ref{lem:KputBall}, for each $l \le k$ and $j=0,1$, there exist $2|I_l|+1$-variate polynomials $\alpha_{j,l}(y,z,u) \in \cK$  such that
\begin{equation}
\label{eq:sparse2}
q_l(x_{I_l}) = \alpha_{0,l}(x_{I_l}+R_le,R_le-x_{I_l},R_l^2-\|x_{I_l}\|^2) + \alpha_{1,l}(x_{I_l}+R_le,R_le-x_{I_l},R_l^2-\|x_{I_l}\|^2)(R_l^2-\|x_{I_l}\|^2).
\end{equation}
The statement then follows by using~\eqref{eq:sparse2} in~\eqref{eq:sparse1}.
\end{proof}

\subsection{The case of SOS polynomials}
\label{sec:SOS}

As exemplified by Schm\"{u}dgen's~\cite{schmudgen1991k} and Putinar's~\cite{putinar1993positive} Positivstellensatz,  as well as the rich, related literature, SOS polynomials form the main polynomial class of choice to construct Positivstellensatz for polynomials~\citep[see also,][for an excellent review of much earlier related work]{reznick2000some}. Not surprisingly, SOS polynomials can be chosen as the class $\cK$ (i.e., setting $\cK = \KSOS$) in the  semi-sparse (resp. fully-sparse) $\cK$-Putinar \germans{} introduced in Theorem~\ref{thm:semiSparse} (resp.  Theorem~\ref{thm:mainSparse}). That is, clearly,
$\KSOMS \subseteq  \KSOS$.
We have then that $\KSOS$ satisfies the assumptions of Theorems~\ref{thm:semiSparse} and~\ref{thm:mainSparse}. Moreover, in this case, the expressions of the semi-sparse (resp. fully-sparse) $\cK$-Putinar \germans{} can be simplified thanks to the
following properties of SOS polynomials. Namely, SOS polynomials are closed under composition

\begin{lemma}\label{lem:simplifySOS}
Let $q\in \R[u_1,\dots,u_s]$, and $h_1,\dots,h_s\in \R[x]$. Then,  if $f\in \KSOS$, then $f(h_1,\dots,h_s)\in \KSOS$.
\end{lemma}

From $\KSOMS \subseteq  \KSOS$, Lemma~\ref{lem:simplifySOS}, and Theorem~\ref{thm:semiSparse}, we can derive a semi-sparse SOS Putinar Positivstellensatz that is basically equivalent to the one introduced in~\citep[][Cor. 7]{kuryatnikova2019copositive}.

\begin{proposition}[Semi-sparse SOS Putinar Positivstellensatz]
\label{prop:KputSOS}
Let $p,g_1,\dots,g_m\in\R[x]$. Then $p(x) > 0$ for all $x \in S$ where $S=\{x \in \B_r: g_j(x) \ge 0, j=1,\dots,m\}$ implies that there exist SOS polynomials $\alpha_l(x)$ for $l=1,2$, and univariate SOS-polynomials $\rho_j(u)$, $j=1,\dots,m$, such that
\begin{equation}
\label{eq:semiSparseSOS}
    p(x) = \alpha_1(x)+\alpha_2(x)(r^2-\|x\|^2) + \sum_{j=1}^m\rho_j(g_j(x))g_j(x).
\end{equation}
\end{proposition}

Note that in obtaining~\eqref{eq:semiSparseSOS} (as well as~\eqref{eq:mainSparseSOS} below), we use that Lemma~\ref{lem:simplifySOS} implies that given the univariate polynomial $\rho \in \KSOS$ and $U \in \R$, it follows that $\rho'(t) := \rho(U-t) \in \KSOS$.

Similarly, from  $\KSOMS \subseteq  \KSOS$, Lemma~\ref{lem:simplifySOS}, and Theorem~\ref{thm:mainSparse}, we can derive a fully-sparse SOS Putinar Positivstellensatz, taking advantage of both term and correlative sparsity in the polynomial and underlying semialgebraic set of interest.

\begin{proposition}[{Fully-sparse SOS Putinar Positivstellensatz~\citep[][Cor. 3.3]{lasserre2006convergent}}]
\label{prop:sparseSOS}
Let $I$ be a cover of $[n]$ satisfying the RIP. Let $\{l_1, l_2, \dots, l_k\} = [k]$.
Let $g_j \in \R[x_{I_{l_j}}]$, $j=1,\dots,m$, and $R \in \R^k_{++}$. Then $p(x) > 0$ on $\{x \in \B^I_R: g_j(x_{I_{l_j}}) \ge 0, j=1,\dots,m\}$ implies that there exist $\SOS$-polynomials $\alpha_{j,l} (x_{I_l})$ for $l=1,\dots,k$ and $j=1,2$, and univariate SOS-polynomials $\sigma_{j}(u)$, for all $j$, such that
\begin{equation}
\label{eq:mainSparseSOS}
p(x) = \sum_{\ell=1}^k[ \alpha_{1,\ell}(x_{I_{\ell}}))+ \alpha_{2,\ell}(x_{I_{\ell}})(r^2-\|x_{I_{\ell}}\|^2)] + \sum_{j=1}^m\sigma_{j}(g_j(x_{I_{l_j}})) g_j(x_{I_{l_j}}).
\end{equation}
\end{proposition}

For any $j=1,\dots,m$ the polynomial $\sigma_j'(x_{I_{l_j}}) = \sigma_{j}(g_j(x_{I_{l_j}})) \in \R[x_{I_{l_j}}]$, belongs to $\KSOS$. From this fact, it follows that Proposition~\ref{prop:sparseSOS} implies the well known fully-sparse SOS Putinar \german{}
derived in~\citep[][Cor. 3.3]{lasserre2006convergent}.

\section{Conclusions and future work}
\label{sec:end}
In this paper, we consider the problem of constructing \germans{} for polynomials over compact semialgebraic sets. Typically, SOS polynomials are used as the base class of polynomials to certify non-negativity when constructing this type of \germans{}. For example, this is the case in the well-known SOS \germans{} of
Schm\"{u}dgen~\cite{schmudgen1991k} and Putinar~\cite{putinar1993positive}. Here, we provide a framework that can be used to generate non-SOS Putinar-type \germans{} for any class of polynomials satisfying a mild assumption; that is, \germans{} resembling the form of Putinar's \german{} in which instead of SOS polynomials, other classes of polynomials are used as the base class to certify non-negativity (see Theorem~\ref{thm:semiSparse} and Theorem~\ref{thm:mainSparse}). Turns out that well-known classes of polynomials such as SONC, SDSOS, DSOS, and \SAGEtext{} polynomials satisfy this assumption. Namely, that the class contains all sum of monomial squares polynomials. Therefore, these classes of polynomials can be used to obtain associated SONC, SDSOS, DSOS, and \SAGEtext{} Putinar-type \german{} (see Corollary~\ref{cor:semiSparse} and Corollary~\ref{cor:mainSparse}).

These results contribute to the literature in polynomial optimization and algebraic geometry in three main ways. First, the results extend results obtained regarding non-SOS \schmu{}-type \germans{}~\citep{dressler2017positivstellensatz, kuryatnikova2019copositive, chandrasekaran2016relative}; that is, \germans{} following the form of \schmu{}'s \german{} in which instead of SOS polynomials, other classes of polynomials are used as the base class to certify non-negativity. By moving from non-SOS \schmu{}-type to non-SOS Putinar-type, one goes from considering \germans{} with an exponential number of terms to \germans{} with a linear number of terms. This is one of the reasons behind the wide popularity of Putinar's \german{} in polynomial optimization, which
is the building block of Lasserre's seminal work in the area~\citep[see, e.g.,][]{anjo12, BlekPT13, lasserre2009moments}. Second, our results show that the new classes of non-SOS polynomials that have been considered in the literature, in part to address the computational expenses associated to checking the membership on the class of SOS polynomials, can be used to construct non-SOS Putinar-type \germans{}. Previous attempts at obtaining this kind of result have failed~\citep{dressler2017positivstellensatz,dressler2018optimization}, in great part because the task is challenging due to the fact that this new classes of non-SOS polynomials lack characteristics of the class of SOS polynomials that are key in deriving SOS \germans{} using algebraic geometry tools~\citep{dressler2017positivstellensatz}. In turn, these results imply that these new non-SOS classes of polynomials can be used to construct approximation hierarchies, that are guaranteed to converge to the solution of polynomial optimization problems, with desired characteristics associated with Putinar-type \germans{}, but without the need to incur in the computational expense associated with checking membership in the class of SOS polynomials.
Thirdly, our results show that the possibility of simplifying \germans{} expression by exploiting sparsity is by no means limited, as in the current literature, to SOS  \germans{}~\citep{lasserre2006convergent, mai2020sparse,wang2020cs,wang2021tssos, kuryatnikova2019copositive}. Instead, our results show that even in the case in which no particular sparsity information is available regarding the polynomials involved in the \germans{}, the obtained non-SOS Putinar-type \germans{} can be written with inherent term sparsity (Theorem~\ref{thm:semiSparse}, Corollary~\ref{cor:semiSparse}, and Remark~\ref{rem:sparsity}). When correlative sparsity characteristics of the polynomials in the \germans{} is known, the  non-SOS Putinar-type \germans{} can be further simplified (Theorem~\ref{thm:mainSparse}, Corollary~\ref{cor:mainSparse}). This takes the study of sparse polynomial optimization and algebraic geometry from SOS centered results, to results that involve a wide class of non-SOS classes of polynomials.

An important aspect of the application of \germans{} are {\em degree bounds}~\citep[see, e.g.,][]{powers2001new}. Given that the size of the corresponding polynomials used in the \germans{} to certify non-negativity
grows exponentially with the number of varaibles and degree of the polynomial whose non-negativity is being certified, having a bound on the required degree of the polynomials constructing the \german{}  may improve the applicability of these \germans{} , and help further close the knowledge gap between $\SOS$ \germans{} and non-SOS \germans{}. This interesting problem is the topic of current future work.

\bibliographystyle{apalike}

\end{document}